\documentclass[12pt,a4]{amsart}
\usepackage{amsmath}
\usepackage{amssymb}
\usepackage{amsthm}
\usepackage{enumitem}
\usepackage{url}
\usepackage{soul}
\usepackage{times}
\usepackage[T1]{fontenc}
\usepackage[export]{adjustbox}

\usepackage{color}
\usepackage{dsfont}
\usepackage{epsfig}
\usepackage[hidelinks]{hyperref}
\usepackage{setspace}
\usepackage{epstopdf}
\usepackage[authoryear,round,sort]{natbib}
\usepackage{comment}
\usepackage{booktabs}
\usepackage{bm}
\usepackage{float}
\usepackage[linesnumbered,ruled,vlined]{algorithm2e}
\SetKwInput{KwInput}{Input}                
\SetKwInput{KwOutput}{Output} 
\usepackage{stmaryrd}
\usepackage{subcaption,graphicx}
\usepackage{float}
\numberwithin{equation}{section}

\usepackage[font=small,labelfont=bf]{caption}
\DeclareCaptionFont{tiny}{\tiny}
\newcommand{\norm}[1]{\left\lVert#1\right\rVert}
\usepackage{geometry}
\usepackage{graphicx}
\theoremstyle{plain}
\newtheorem{proposition}{Proposition}[section]

\newtheorem{theorem}{Theorem}[section]
\newtheorem{lemma}{Lemma}[section]
\newtheorem{corollary}{Corollary}[section]

\theoremstyle{definition}
\newtheorem{definition}{Definition}[section]

\newtheorem{example}{Example}[section]

\theoremstyle{remark}
\newtheorem{rk}{Remark}[section]
\expandafter\let\expandafter\oldproof\csname\string\proof\endcsname
\let\oldendproof\endproof
\renewenvironment{proof}[1][\proofname]{%
  \oldproof[\noindent\textbf{#1.} ]%
}{\oldendproof}

\newcommand{\be}{\begin{equation}}
\newcommand{\ee}{\end{equation}}
\newcommand{\by}{\begin{eqnarray*}}
\newcommand{\ey}{\end{eqnarray*}}
\DeclareMathOperator*{\argmax}{arg\,max}
\DeclareMathOperator*{\argmin}{arg\,min}
\DeclareMathOperator*{\argsup}{arg\,sup}
\DeclareMathOperator*{\arginf}{arg\,inf}

\renewcommand{\leq}{\leqslant}
\renewcommand{\geq}{\geqslant}
\usepackage{xcolor}
\definecolor{dark-red}{rgb}{0.4,0.15,0.15}
\definecolor{dark-blue}{rgb}{0.15,0.15,0.4}
\definecolor{medium-blue}{rgb}{0,0,0.5}
\hypersetup{
    colorlinks, linkcolor={red},
    citecolor={blue}, urlcolor={blue}
}
\allowdisplaybreaks
\title[Markov chain entropy games and the geometry of their Nash equilibria]{Markov chain entropy games and the geometry of their Nash equilibria}
\author{Michael C.H. Choi}
\address{Department of Statistics and Data Science, National University of Singapore, Singapore}
\email{mchchoi@nus.edu.sg}

\author{Geoffrey Wolfer}
\address{Waseda University Center for Data Science, Tokyo, Japan}
\email{geo.wolfer@aoni.waseda.jp}

\begin{document}

\date{\today}
\maketitle

\begin{abstract}
	\textcolor{black}{We introduce and study a two-player zero-sum game between a probabilist and Nature defined by a convex function \( f \), a finite collection \( \mathcal{B} \) of Markov generators (or its convex hull), and a target distribution \( \pi \). The probabilist selects a mixed strategy \( \mu \in \mathcal{P}(\mathcal{B}) \), the set of probability measures on $\mathcal{B}$, while Nature adopts a pure strategy and selects a \( \pi \)-reversible Markov generator \( M \). The probabilist receives a payoff equal to the \( f \)-divergence \( D_f(M \| L) \), where \( L \) is drawn according to \( \mu \). We prove that this game always admits a mixed strategy Nash equilibrium and satisfies a minimax identity. In contrast, a pure strategy equilibrium may fail to exist. We develop a projected subgradient method to compute approximate mixed strategy equilibria with provable convergence guarantees. Connections to information centroids, Chebyshev centers, and Bayes risk are discussed. This paper extends earlier minimax results on \( f \)-divergences such as \citep{H97,HO97,GZ06} to the context of Markov generators.}
	
	\smallskip
	
	
\end{abstract}

\tableofcontents


\section{Introduction}\label{sec:intro}

\textcolor{black}{This paper studies a two-player zero-sum game between a probabilist and Nature, defined by parameters \((f, \mathcal{B}, \pi)\). These will be introduced formally in Sections \ref{sec:prelim} and \ref{sec:main}, but we outline them here to motivate the setup. The function \(f\) is convex and satisfies standard regularity conditions; \(\mathcal{B}\) is either a finite set \(\{L_i\}_{i=1}^n\) or its convex hull, with each \(L_i\) a Markov generator on a finite state space \(\mathcal{X}\); and \(\pi\) is a given target distribution on \(\mathcal{X}\). We analyze two versions of the game: one with pure strategies, and one with mixed strategies. In the mixed strategy setting, the probabilist selects a probability measure \(\mu \in \mathcal{P}(\mathcal{B})\) and samples \(L \in \mathcal{B}\) according to \(\mu\), while Nature adopts a pure strategy and chooses a reversible Markov generator \(M \in \mathcal{L}(\pi)\). The probabilist receives a payoff of \(D_f(M \| L)\), the \(f\)-divergence from \(L\) to \(M\). In the pure strategy variant, both players choose deterministically from their respective strategy sets.}


\textcolor{black}{Games of a similar flavor appear in the literature under the name ``statistician versus Nature'' problems \citep{H97,HO97,GZ06}, where the statistician pays a loss based on divergence between the estimator and the ground truth distribution chosen by Nature. Our game reverses this direction: the probabilist receives a gain from Nature. This role reversal is motivated by differences between reversible and non-reversible Markov chains as follows:}


\textcolor{black}{Reversibility is fundamental in many physical systems---quantum and classical mechanics are both time-reversible \citep{S1931}. Accordingly, we model Nature as selecting from the set of \(\pi\)-reversible generators. This also reflects common practice in modern Markov chain Monte Carlo methods, where algorithms such as Metropolis-Hastings and overdamped Langevin diffusions are \(\pi\)-reversible.}


\textcolor{black}{In contrast, non-reversible chains have attracted recent interest for their potential to accelerate convergence. Algorithms of this type appear in various works (e.g. \cite{Hwang93,Hwang05,RR15,Bie16,DHN00,BS16,G98,KS23}). We therefore let the probabilist choose from a possibly non-reversible inventory of Markov generators \(\mathcal{B}\). The quantity \(D_f(M \| L)\) can be interpreted as an information-theoretic advantage or performance gain due to non-reversibility.}


While the above explanations justify the role of the probabilist and Nature in the entropy game, from a mathematical point of view however, these two roles can be safely interchanged, or one may wish to replace the term "probabilist" by "Player A" and "Nature" by "Player B" throughout the entire manuscript.

\textcolor{black}{This setting raises several natural questions: Does a Nash equilibrium exist \citep{KarlinPeres17,MSZ20}? Is it unique? Can it be computed efficiently? What is the game’s value? Is there a sequential analogue?}

\textcolor{black}{While we do not aim to answer all of these questions, this paper focuses on several foundational aspects: the existence and uniqueness of an equilibrium, cases where equilibria may fail to exist, and an optimization algorithm to approximate them efficiently. We also explore connections with the Bayes risk and geometric notions such as information centroids and Chebyshev centers.}

\textcolor{black}{Further motivation comes from a question raised by Laurent Miclo about whether the Metropolis-type generators
\[
P_{-\infty}(x,y) := \min\{L(x,y),L_{\pi}(x,y)\}, \quad
P_{\infty}(x,y) := \max\{L(x,y),L_{\pi}(x,y)\}, \quad x \ne y,
\]
can be interpreted game-theoretically. Here, \(L_{\pi}\) is the \(\pi\)-dual of \(L\), and diagonal terms are set so that each row sums to zero (see Section~\ref{sec:prelim}). While \(P_{-\infty}\) corresponds to the classical continuous-time Metropolis-Hastings generator, \(P_{\infty}\) has been studied in \cite{DM09,Choi16,CH18}. Since both involve pointwise minimum or maximum operations, a game-theoretic interpretation is natural. As we show in Example~\ref{ex:1mix}, when \(D_f\) is taken to be the total variation distance, both \(P_{-\infty}\), \(P_{\infty}\), and their convex combinations can be understood as part of a mixed strategy Nash equilibrium in the entropy game. This gives a possible answer to Miclo’s question.}

We summarize our main contributions as \textcolor{black}{follows}:
\begin{enumerate}
	\item \textcolor{black}{We introduce} the two-person Markov chain entropy games of \textcolor{black}{the} probabilist against Nature. This naturally generalizes the game of \textcolor{black}{the} statistican against Nature \citep{H97,HO97,GZ06}, where instead of Markov generators, the game in these references involves probability measures. 
	
	\textcolor{black}{We draw a perhaps surprising connection between reversiblizations of Markov generators, widely studied in modern theory \citep{M97,Fill91,Paulin15,Choi16}, game-theoretic concepts, and geometry.}
%
%
%
%
%
	
	\item \textcolor{black}{We introduce} the notion of a weighted information centroid in the context of Markov generators and \textcolor{black}{highlight its central role in analyzing Nash equilibria} of entropy games. This naturally generalizes the notion of \textcolor{black}{an} information centroid of Markov chains introduced by \cite{CW23}. Our analysis shows that the (mixed strategy) Nash equilibrium is intimately related to the notion of Chebyshev center, which can be interpreted as a specific weighted information centroid. This important observation allows us to analyze the existence and uniqueness of a Nash equilibrium in the game.
	
	
	\item \textcolor{black}{We propose} and analyze a simple projected subgradient algorithm to find an approximate Nash equilibrium with a provable convergence rate of $\mathcal{O}(1/\sqrt{t})$. A central question in game theory, in particular algorithmic game theory, lies in developing efficient algorithms for (approximate) computation of the Nash equilibrium. To this end, we propose a simple and easy-to-implement projected subgradient algorithm which \textcolor{black}{takes advantage of} the information geometry of the underlying Markov generators.
	
\end{enumerate}
The rest of this paper is organized as \textcolor{black}{follows}. In Section \ref{sec:prelim}, \textcolor{black}{we introduce and recall various notions}, in particular the notion of weighted information centroid, which plays a central role in our subsequent investigation. We proceed to present the main results of the paper, where we first address fundamental questions concerning Chebyshev \textcolor{black}{centers}, weighted information centroids and minimax values in Section \ref{sec:main}, followed by a pure strategy game-theoretic analysis in Section \ref{subsec:puregame}. The corresponding setting of mixed strategy game is presented in Section \ref{subsec:mixgame}, and several simple yet illustrative examples are given in Section \ref{subsubex:pureex} and Section \ref{subsubex:mixedex}. The design and analysis of a novel projected subgradient algorithm to find an approximate mixed strategy Nash equilibrium is stated in Section \ref{subsec:algo}. Finally, the proofs of the main results are \textcolor{black}{presented} in Section \ref{sec:proofsmain}.

\subsection{Notations}

In this subsection we introduce some commonly used notations throughout the manuscript. For $a,b \in \mathbb{Z}$, we write $\llbracket a,b\rrbracket := \{a,a+1,\ldots,b-1,b\}$ and $\llbracket n \rrbracket := \llbracket 1,n \rrbracket$ for $n \in \mathbb{N}$. We denote by $\mathbf{0}$ the all-zeros matrix on the finite state space $\mathcal{X}$. For a given function \(g(n)\), we say that it is \(\mathcal{O}(h(n))\) if there exist constants \(C > 0\) and \(n_0\) such that \(g(n) \leq C h(n)\) for all \(n \geq n_0\). 

\section{Preliminaries}\label{sec:prelim}

\subsection{Weighted information centroids}\label{subsec:weightedic}

Consider a convex function \(f : \mathbb{R}_+ \rightarrow \mathbb{R}_+\) satisfying \(f(1) = 0\), and a given positive discrete distribution $\pi$ on a finite state space \(\mathcal{X}\), \textcolor{black}{where $\mathbb{R}_+$ is the set of non-negative real numbers}. Define \(\mathcal{L}\) as the set of Markov infinitesimal generators on \(\mathcal{X}\). These generators correspond to \(\mathcal{X} \times \mathcal{X}\) matrices having non-negative off-diagonal elements and zero row sums. A generator \(L\) is said to be \(\pi\)-stationary if \(\pi L = 0\). Moreover, \(\mathcal{L}(\pi) \subset \mathcal{L}\) is the subset of \(\pi\)-reversible generators, where we recall that a generator $L$ is said to be $\pi$-reversible if it satisfies the detailed balance condition, that is, for all $x \neq y \in \mathcal{X}$, we have
$\pi(x) L(x,y) = \pi(y) L(y,x).$

In view of \cite[Proposition 1.2]{JK14}, we define the \(\pi\)-dual of a generator \(L \in \mathcal{L}\) to be \(L_{\pi}\). For \(x \neq y\), the off-diagonal elements of $L_{\pi}$ are given by
\[
L_{\pi}(x,y) = \frac{\pi(y)}{\pi(x)}L(y,x),
\]
and the diagonal ones ensure zero row sums for all rows. In the special case when \(L\) has \(\pi\) as its unique stationary distribution, we then have \(L_{\pi} = L^*\), the adjoint of \(L\) in \(\ell^2(\pi)\) or its time-reversal. Here the space \(\ell^2(\pi)\) is the standard weighted \(\ell^2\) Hilbert space endowed with the inner product \(\langle \cdot, \cdot \rangle_{\pi}\) given by, for any functions $g,h: \mathcal{X}\rightarrow\mathbb{R}$,
\begin{align*}
	\langle g,h\rangle_\pi:=\sum_{x\in \mathcal{X}}g(x)h(x)\pi(x).
\end{align*}

Following \cite{DM09}, for a given target \(\pi\) and Markov infinitesimal generators \(M, L \in \mathcal{L}\), the \(f\)-divergence from $L$ to $M$ with respect to $\pi$ is defined as
\begin{equation}\label{def:fdivML}
	D_f(M || L) =  \sum_{x \in \mathcal{X}} \pi(x) \sum_{y \in \mathcal{X}\setminus\{x\}} L(x,y) f\left(\frac{M(x,y)}{L(x,y)}\right),
\end{equation}
where the convention of \(0 f(a/0) := 0\) for \(a \geq 0\) applies. If \(f^*\) is the convex conjugate of \(f\), defined by \(f^*(t) = tf(1/t)\) for \(t > 0\), then
\[
D_f(M || L) = D_{f^*}(L || M),
\]
and \(f^*(1) = 0\). \textcolor{black}{In particular}, when \(f\) is self-conjugate, that is \(f^* = f\), the \(f\)-divergence in \eqref{def:fdivML} is symmetric.

For a general generator \(L\) not necessarily with \(\pi\) as its stationary distribution, we consider projecting \(L\) onto \(\mathcal{L}(\pi)\) under \(f\)-divergence \(D_f\). Following \cite{WW21,CW23}, the notions of \(f\)-projection and \(f^*\)-projection are defined as:
\begin{equation}\label{def:emprojection}
	M^f = M^f(L,\pi) := \argmin_{M \in \mathcal{L}(\pi)} D_f(M || L), \quad M^{f^*} = M^{f^*}(L,\pi) = \argmin_{M \in \mathcal{L}(\pi)} D_f(L || M).
\end{equation}

To \textcolor{black}{state} our main results, we now introduce the notion of \textcolor{black}{a weighted information centroid} of Markov chains. It generalizes the notion \textcolor{black}{of an (unweighted) information centroid} of Markov chains \textcolor{black}{introduced} in \cite{CW23}. Let $\mathbf{w} = (w_1,\ldots,w_n) \in \mathbb{R}^n_+$ be a weight vector in the probability simplex $\mathcal{S}_n$ of $n$ elements, that is, 
\begin{align}\label{def:probsim}
	\mathcal{S}_n := \bigg\{\mathbf{w} = (w_1,\ldots,w_n) \in \mathbb{R}^n_+; ~\sum_{i=1}^n w_i = 1\bigg\}.
\end{align}
Note that we denote the simplex by $\mathcal{S}_n$ instead of $\mathcal{P}(\llbracket n \rrbracket)$ throughout the manuscript. Given a sequence of Markov generators $\{L_i\}_{i=1}^n$, where $L_i \in \mathcal{L}$ for $i \in \llbracket n \rrbracket$, we define the notions of $\mathbf{w}$-weighted $f$-projection centroid and $f^*$-projection centroid to be respectively
\begin{align*}
	M^{f}_n &= M^{f}_n(\mathbf{w},\{L_i\}_{i=1}^n,\pi) := \argmin_{M \in \mathcal{L}(\pi)} \sum_{i=1}^n w_i D_f(M || L_i), \\
	M^{f^*}_n &= M^{f^*}_n(\mathbf{w},\{L_i\}_{i=1}^n,\pi) = \argmin_{M \in \mathcal{L}(\pi)} \sum_{i=1}^n w_i D_f(L_i || M).
\end{align*}
We consider two important special cases: first in the case with $n = 1$, the above notions reduce to $M_{1}^f = M^f$ and $M_{1}^{f^*} = M^{f^*}$ respectively as introduced in \eqref{def:emprojection}. In the second special case, let $n \in \mathbb{N}$ and $\{\mathbf{e}_i\}_{i=1}^n$ be the standard unit vectors. It is \textcolor{black}{immediate} that $M^f_n(\mathbf{e}_i,\{L_i\}_{i=1}^n,\pi) = M^f(L_i,\pi)$ and $M^{f^*}_n(\mathbf{e}_i,\{L_i\}_{i=1}^n,\pi) = M^{f^*}(L_i,\pi)$. 

Our first main result in this Section establishes existence and uniqueness of $f$ and $f^*$-projection centroids under strict convexity of $f$, and its proof is given in Section \ref{subsubsec:pfexistunique}.

\begin{theorem}[Existence and uniqueness of weighted $f$ and $f^*$-projection centroids]\label{thm:existuniquecentroid}
	Suppose we are given a sequence of Markov generators $\{L_i\}_{i=1}^n$, where $L_i \in \mathcal{L}$ for $i \in \llbracket n \rrbracket$, a $f$-divergence $D_f$ generated by a strictly convex $f$ which is assumed to have a derivative at $1$ given by $f^{\prime}(1) = 0$, and a weight vector $\mathbf{w} \in \mathcal{S}_n$. We further assume that there exists at least one $L_i \neq \mathbf{0}$ and at least one $w_i > 0$ when $L_i \neq \mathbf{0}$. A weighted $f$-projection of $D_f$ (resp.~$f^*$-projection of $D_{f^*}$) that minimizes the mapping 
	$$\mathcal{L}(\pi) \ni M \mapsto \sum_{i=1}^n w_i D_f(M||L_i) \quad \left(\textrm{resp.}~=\sum_{i=1}^n w_i D_{f^*}(L_i||M)\right)$$
	\textcolor{black}{uniquely exists and is denoted by} $M^f_{n}$. 
\end{theorem}

\begin{rk}
	By replacing $f$ by $f^*$ in the main result of Theorem \ref{thm:existuniquecentroid}, we can analogously define $M^{f^*}_n$. Precisely, a weighted $f^*$-projection of $D_f$ (resp.~$f$-projection of $D_{f^*}$) that minimizes the mapping 
	$$\mathcal{L}(\pi) \ni M \mapsto \sum_{i=1}^n w_i D_f(L_i||M) \quad \left(\textrm{resp.}~=\sum_{i=1}^n w_iD_{f^*}(M||L_i)\right)$$
	\textcolor{black}{uniquely exists and is denoted by} $M^{f^*}_{n}$.
\end{rk}

\begin{rk}[On the convexity of $f$]
	We emphasize that throughout this manuscript, unless otherwise specified, $f$ is assumed to be a convex function rather than a strictly convex function. This yields the following interesting consequences. Suppose that $f$ is a convex function with $f(s) = 0$ for some $s > 0, s \neq 1$. For a given $L \in \mathcal{L}$, we thus have
	$$D_f(sL||L) = 0,$$
	and hence the two generators $sL$ and $L$ are indistinguishable with respect to $D_f$. From a transition semigroup perspective, this is justifiable since $P^t := e^{sLt}$, the transition semigroup generated by $sL$, is merely a time-change of $Q^t := e^{Lt}$, the transition semigroup generated by $L$, for $t \geq 0$.
\end{rk}
\begin{rk}
	We discuss the importance of the additional assumptions in Theorem \ref{thm:existuniquecentroid} \textcolor{black}{here}.
	In the first case, if $L_i = \mathbf{0}$ for all $i \in \llbracket n \rrbracket$, then obviously for any $M \in \mathcal{L}(\pi)$ we have $D_f(M||L_i) = 0$ and hence
	$$\sum_{i=1}^n w_i D_f(M||L_i) = 0.$$
	As such, this shows the existence of weighted $f$-projection centroids \textcolor{black}{which} are not unique even if $f$ is strictly convex.
	In the second case, suppose that there exists at least one $L_i \neq \mathbf{0}$ and $w_i = 0$ whenever $L_i \neq \mathbf{0}$. This again leads to, for any $M \in \mathcal{L}(\pi)$,
	$$\sum_{i=1}^n w_i D_f(M||L_i) = 0.$$
	We thus exclude the above two degenerate cases in Theorem \ref{thm:existuniquecentroid}.
\end{rk}

\begin{rk}
	We stress that in Theorem \ref{thm:existuniquecentroid}, the given positive distribution $\pi$ is arbitrary and at the same time the generators $L_i$ need not admit $\pi$ as a stationary distribution or even be irreducible in the first place. The purpose of Theorem \ref{thm:existuniquecentroid} is to consider the joint projection of these generators $L_i$ onto the space $\mathcal{L}(\pi)$.
\end{rk}
In the second result of this Section, we explicitly calculate the weighted $f$ and $f^*$-projection centroids $M^f_{n}$ and $M^{f^*}_{n}$ under some common $f$-divergences. The proof is deferred to Section \ref{subsubsec:examplecentroid}. For these common choices of $f$, as shown in \cite{CW23}, the power mean reversibilizations $P_p$ appear naturally as the corresponding $f$ or $f^*$-projection. For \(x \neq y \in \mathcal{X}\) and \(p \in \mathbb{R}\setminus\{0\}\), the off-diagonal entries of $P_p$ are given by
\begin{equation}\label{def:Pp}
	P_{p}(x,y) := \left(\frac{L(x,y)^{p} + L_{\pi}(x,y)^{p}}{2}\right)^{1/p},
\end{equation}
where diagonal entries ensure zero row sums. The limiting cases are defined to be
\begin{align}
	P_0(x,y) &:= \lim_{p \to 0} P_{p}(x,y) = \sqrt{L(x,y)L_{\pi}(x,y)}, \nonumber \\
	P_{\infty}(x,y) &:= \lim_{p \to \infty} P_{p}(x,y) = \max\{L(x,y),L_{\pi}(x,y)\}, \label{eq:Pinfty}\\
	P_{-\infty}(x,y) &:= \lim_{p \to -\infty} P_{p}(x,y) = \min\{L(x,y),L_{\pi}(x,y)\}. \label{eq:P-infty}
\end{align}
Note that in the special case when the weights are given by $w_i = 1/n$ for all $i \in \llbracket n \rrbracket$, we recover the results in \cite{CW23}, as expected.

\begin{theorem}[Examples of weighted $f$ and $f^*$-projection centroids]\label{thm:examplecentroid}
	Given a sequence of Markov generators $(L_i)_{i=1}^n$, where $L_i \in \mathcal{L}$ for $i \in \llbracket n \rrbracket$ and a weight vector $\mathbf{w} \in \mathcal{S}_n$. We further assume that there exists at least one $L_i \neq \mathbf{0}$ and at least one $w_i > 0$ when $L_i \neq \mathbf{0}$. Recall the power mean reversiblizations $P_p$ in \eqref{def:Pp}.
	\begin{enumerate}
		\item(weighted $f$ and $f^*$-projection centroids under $\alpha$-divergence)\label{it:alphacentroid}
		Let $f(t) = \frac{t^{\alpha} - \alpha t - (1-\alpha)}{\alpha(\alpha-1)}$ for $\alpha \in \mathbb{R}\backslash\{0,1\}$. The unique $f$-projection centroid $M^{f}_n$ is given by, for $x \neq y \in \mathcal{X}$,
		\begin{align*}
			M^{f}_n(x,y) =\left(\sum_{i=1}^n w_i \left( M^{f}(L_i,\pi)(x,y)\right)^{1-\alpha}\right)^{1/(1-\alpha)},
		\end{align*}
		while the unique $f^*$-projection centroid $M^{f^*}_n$ is given by, for $x \neq y \in \mathcal{X}$,
		\begin{align*}
			M^{f^*}_n(x,y) =\left(\sum_{i=1}^n w_i\left(M^{f^*}(L_i,\pi)(x,y)\right)^{\alpha}\right)^{1/\alpha},
		\end{align*}
		where $M^{f}, M^{f^*}$ are respectively the $P_{1-\alpha}, P_{\alpha}$-reversiblization.
		
%
		\item(weighted $f$ and $f^*$-projection centroids under squared Hellinger distance)\label{it:hellingercentroid}
		Let $f(t) = (\sqrt{t}-1)^2$.
		The unique $f$-projection centroid $M^{f}_n$ is given by, for $x \neq y \in \mathcal{X}$,
		\begin{align}\label{eq:cephellinger}
			M^{f}_n(x,y) =\left(\sum_{i=1}^n w_i\sqrt{M^{f}(L_i,\pi)(x,y)}\right)^2,
		\end{align}
		while the unique $f^*$-projection centroid $M^{f^*}_n$ is given by, for $x \neq y \in \mathcal{X}$,
		\begin{align*}
			M^{f^*}_n(x,y) =\left(\sum_{i=1}^n w_i\sqrt{M^{f^*}(L_i,\pi)(x,y)}\right)^2,
		\end{align*}
		where $M^{f^*} = M^f$ is the $P_{1/2}$-reversiblization.
	
		\item($f$ and $f^*$-projection centroids under Kullback-Leibler divergence)
		Let $f(t) = t \ln t - t + 1$. The unique $f$-projection centroid $M^{f}_n$ is given by, for $x \neq y \in \mathcal{X}$,
		\begin{align*}
			M^{f}_n(x,y) = \prod_{i=1}^n \left(M^{f}(L_i,\pi)(x,y)\right)^{w_i},
		\end{align*}
		where $0^0 := 0$ in the expression above, while the unique $f^*$-projection centroid $M^{f^*}_n$ is given by, for $x \neq y \in \mathcal{X}$,
		\begin{align*}
			M^{f^*}_n(x,y) = \sum_{i=1}^n w_i M^{f^*}(L_i,\pi)(x,y),
		\end{align*}
		where $M^{f}, M^{f^*}$ are respectively the $P_{0}, P_{1}$-reversiblization.
	\end{enumerate}
\end{theorem}

\begin{rk}
	Item \eqref{it:hellingercentroid} can be considered as a special case of item \eqref{it:alphacentroid} by taking $\alpha = 1/2$. Another important special case of $\alpha$-divergence is the $\chi^2$-divergence where we take $\alpha = 2$.
\end{rk}

\subsection{Two-person pure strategy and mixed strategy games}\label{subsec:introgame}

\textcolor{black}{This subsection reviews some classical game-theoretic notions in the context of the probabilist against Nature game. Much of the exposition in this subsection are drawn from content in \cite[Section $5.2.5$ and $5.4.3$]{Boyd2004} and \cite[Chapter $2$]{KarlinPeres17}. Consider the following two-person non-cooperative zero-sum pure strategy game with respect to the parameters $(f, \mathcal{B},\pi)$. If Nature chooses $M \in \mathcal{L}(\pi)$, the \emph{strategy set of Nature}, while the probabilist chooses $L \in \mathcal{B}$, the \emph{strategy set of the probabilist}, then Nature pays an amount $D_f(M||L)$ to the probabilist. As a result, Nature wants to minimize $D_f$ while the probabilist seeks to maximize $D_f$. The quantity $D_f(M||L)$ is known as the \emph{payoff function} of the game, and it is said to be a \emph{pure strategy} game as both players are only allowed to choose deterministically from their respective strategy sets. On the other hand, a game is said to be a \emph{mixed strategy} game if the players are allowed to choose randomly within their respective strategy sets. In the sequel, we shall consider $\mathcal{B}$, the strategy set of the probabilist, to be either $\{L_i\}_{i=1}^n$ or its convex hull $\mathrm{conv}((L_i)_{i=1}^n)$. From now on we shall call this a two-person pure strategy game.}

Suppose that Nature chooses its strategy $M$ first, followed by the probabilist, who has the knowledge of the choice of Nature. The probabilist thus seeks to choose $L \in \mathcal{B}$ to maximize the payoff $D_f(M||L)$. The resulting payoff is $\sup_{L \in \mathcal{B}} D_f(M||L)$, which depends on $M$, the choice of Nature. Nature assumes that the probabilist will choose this strategy, and hence Nature seeks to minimize the worst-case payoff, that is, Nature chooses a strategy in the set
$$\argmin_{M \in \mathcal{L}(\pi)} \sup_{L \in \mathcal{B}} D_f(M||L),$$
known as a \emph{minimax strategy of Nature}, which yields
$$\overline{v}^f = \overline{v}^f(\mathcal{B},\pi) := \inf_{M \in \mathcal{L}(\pi)} \sup_{L \in \mathcal{B}} D_f(M||L),$$
the payoff from Nature to the probabilist. $\overline{v}^f$ is also known as the \emph{minimax value} of the game.

Now, suppose that the order of play is reversed: the probabilist chooses $L \in \mathcal{B}$ first, and Nature, with the knowledge of $L$, picks from the strategy set $\mathcal{L}(\pi)$. Following a similar argument as before, if both players follow the optimal strategy, the probabilist chooses a strategy in the set 
$$\argsup_{L \in \mathcal{B}} \inf_{M \in \mathcal{L}(\pi)} D_f(M||L),$$
known as a \emph{maximin strategy of the probabilist}, which yields
$$\underline{v}^f = \underline{v}^f(\mathcal{B},\pi) :=  \sup_{L \in \mathcal{B}} \inf_{M \in \mathcal{L}(\pi)} D_f(M||L),$$
the payoff from Nature to the probabilist. $\underline{v}^f$ is known as the \emph{maximin value} of the game.

\textcolor{black}{If the minimax equality holds, that is, the minimax and maximin value coincides with
$$\overline{v}^f = \underline{v}^f =: v = v(f,\mathcal{B},\pi),$$
then $v$ is the \emph{value of the pure strategy game}. Let $M^f$ and $L^f$ be respectively a minimax strategy of Nature and a maximin strategy of the probabilist. The pair $(M^f,L^f)$ is called a \emph{pure strategy Nash equilibrium} with respect to $(f,\mathcal{B},\pi)$. $M^f$ is said to be an optimal strategy of Nature while $L^f$ is said to be an optimal strategy of the probabilist.}

\textcolor{black}{The above pure strategy game can be generalized to a situation where the probabilist uses a mixed strategy while Nature maintains a pure strategy. From now on we shall refer this to a two-person mixed strategy game with respect to the parameters $(f, \mathcal{B},\pi)$. In this mixed strategy game, the probabilist first chooses a \emph{prior measure} $\mu \in \mathcal{P}(\mathcal{B})$ and picks a $L \in \mathcal{B}$ at random according to $\mu$. On the other hand Nature, knowing $\mu$ but not knowing $L$, follows a pure strategy to choose $M \in \mathcal{L}(\pi)$. Afterwards, Nature pays an amount $D_f(M||L)$ to the probabilist. This serves as a natural generalization of the game of a statistican against Nature proposed and developed in \cite{H97,HO97,GZ06}.}

With these notions in mind, we now introduce the minimax and maximin value of this game:

\begin{definition}[Minimax and maximin values in the two-person mixed strategy game]\label{def:minimaxmix}
	Consider the two-person mixed strategy game. The \emph{minimax value} and the \emph{maximin value} of this game are defined respectively to be
	\begin{align}
		\overline{V}^f = \overline{V}^f(\mathcal{B},\pi) &:= \inf_{M \in \mathcal{L}(\pi)} \sup_{\mu \in \mathcal{P}(\mathcal{B})} \int_{\mathcal{B}} D_f(M||L)\, \mu(dL), \label{eq:overlineVf}\\
		\underline{V}^f = \underline{V}^f(\mathcal{B},\pi) &:= \sup_{\mu \in \mathcal{P}(\mathcal{B})} \inf_{M \in \mathcal{L}(\pi)} \int_{\mathcal{B}} D_f(M||L)\, \mu(dL). \label{eq:underlineVf}
	\end{align}
	Note that $\overline{V}^f \geq \underline{V}^f$. We can define analogously $\overline{V}^{f^*}, \underline{V}^{f^*}$ by replacing $f$ by $f^*$ above.
\end{definition}

Next, we define the notions of Bayes risk and Bayes strategy of Nature:

\begin{definition}[Bayes risk and Bayes strategy]
	Given a $\mu \in \mathcal{P}(\mathcal{B})$, we say that a Markov generator $M_{\mu} \in \mathcal{L}(\pi)$ is a \emph{Bayes strategy} with respect to $\mu$ if the mapping $M \mapsto  \int_{\mathcal{B}} D_f(M||L)\, \mu(dL)$ attains its infimum at $M_{\mu}$, that is,
	$$ \inf_{M \in \mathcal{L}(\pi)} \int_{\mathcal{B}} D_f(M||L)\, \mu(dL) = \int_{\mathcal{B}} D_f(M_{\mu}||L)\, \mu(dL).$$
	The minimum value $ \int_{\mathcal{B}} D_f(M_{\mu}||L)\, \mu(dL)$ is said to be the \emph{Bayes risk} with respect to $\mu$,
\end{definition}

Finally, we give the notions of minimax and maximin strategy of this game:

\begin{definition}[Minimax and maximin strategies in the two-person mixed strategy game]\label{def:minimaxmixstrategy}
	Consider the two-person mixed strategy game. A strategy in the set
	$$\arginf_{M \in \mathcal{L}(\pi)}  \sup_{\mu \in \mathcal{P}(\mathcal{B})} \int_{\mathcal{B}} D_f(M||L)\, \mu(dL)$$
	is known as a \emph{minimax strategy of Nature}, while a strategy in the set
	$$\argsup_{\mu \in \mathcal{P}(\mathcal{B})} \inf_{M \in \mathcal{L}(\pi)} \int_{\mathcal{B}} D_f(M||L)\, \mu(dL)$$
	is known as a \emph{maximin strategy of the probabilist}.
\end{definition}

\textcolor{black}{If the minimax equality holds, that is, 
$$\overline{V}^f(\mathcal{B},\pi) = \underline{V}^f(\mathcal{B},\pi) =: V = V(f, \mathcal{B},\pi),$$
then we say that $V$ is the \emph{value of the mixed strategy game}. The pair of minimax strategy of Nature and maximin strategy of the probabilist is a \emph{mixed strategy Nash equilibrium} with respect to the parameters $(f, \mathcal{B},\pi)$.}

\section{Main results}\label{sec:main}

In this Section, we gradually introduce and state the main results of this paper. We first state two well-known and useful max-min inequalities for $D_f$ and $D_{f^*}$. Let $\mathcal{A}, \mathcal{B} \subseteq \mathcal{L}$. We have
\begin{align}
	\sup_{L \in \mathcal{B}} \inf_{M \in \mathcal{A}} D_f(M||L) &\leq \inf_{M \in \mathcal{A}} \sup_{L \in \mathcal{B}}  D_f(M||L), \label{eq:minimaxineq1}\\
	\sup_{L \in \mathcal{B}} \inf_{M \in \mathcal{A}} D_f(L||M) &\leq \inf_{M \in \mathcal{A}} \sup_{L \in \mathcal{B}}  D_f(L||M). \label{eq:minimaxineq2}
\end{align}


In the following, we are primarily interested in projection onto the space $\mathcal{L}(\pi)$, and hence we shall take $\mathcal{A} = \mathcal{L}(\pi)$. \textcolor{black}{We now illustrate the relationship between various minimax and maximin values and saddle point}:

\begin{definition}[Minimax and maximin values, saddle point property and saddle point]\label{def:minimaxvalues}
	Let $\mathcal{B} \subseteq \mathcal{L}$. \textcolor{black}{Recall that}
	\begin{align*}
		\overline{v}^f = \overline{v}^f(\mathcal{B},\pi) &= \inf_{M \in \mathcal{L}(\pi)} \sup_{L \in \mathcal{B}} D_f(M||L), \\
		\underline{v}^f = \underline{v}^f(\mathcal{B},\pi) &= \sup_{L \in \mathcal{B}} \inf_{M \in \mathcal{L}(\pi)} D_f(M||L).
	\end{align*}
	We say that $\overline{v}^f$ is the minimax value with respect to $(f,\mathcal{B},\pi)$, and similarly $\underline{v}^f$ is the maximin value with respect to $(f,\mathcal{B},\pi)$. In view of \eqref{eq:minimaxineq1} and \eqref{eq:minimaxineq2}, we have
	\begin{align*}
		\overline{v}^f \geq \underline{v}^f, \quad \overline{v}^{f^*} \geq \underline{v}^{f^*}.
	\end{align*}
	If the equality holds, that is,
	$$\overline{v}^f = \underline{v}^f,$$
	we say that the saddle point property is satisfied with respect to $(f,\mathcal{B},\pi)$. In this case, suppose that the pair $(M^f,L^f) \in \mathcal{L}(\pi) \times \mathcal{L}$ attains the optimal value, that is,
	$$\overline{v}^f = D_f(M^f||L^f),$$
	then we say that $(M^f,L^f)$ is a saddle point with respect to $(f,\mathcal{B},\pi)$.
\end{definition}

\textcolor{black}{Let us now assume} that we are given a set of Markov generators $\{L_i\}_{i=1}^n$, where $L_i \in \mathcal{L}$ for $i \in \llbracket n \rrbracket$, and in the sequel we shall investigate either $\mathcal{B} = \{L_i\}_{i=1}^n$ or the convex hull of $\{L_i\}_{i=1}^n$. Consider the minimax problem 
\begin{align}\label{problem:minimax}
	\inf_{M \in \mathcal{L}(\pi)} \max_{i \in \llbracket n \rrbracket} D_f(M||L_i).
\end{align}
As $L \mapsto D_f(M||L)$ is convex for fixed $M$ and pointwise maximum of convex functions preserves convexity \cite[Section $3.2.3$]{Boyd2004}, the mapping $M \mapsto \max_{i \in \llbracket n \rrbracket} D_f(M||L_i)$ is thus convex. As such the outer minimization is a convex minimization problem over the convex set $\mathcal{L}(\pi)$. Inspired by the reformulation of an analogous minimax problem in the context of probability measures in \cite[equation $(3)$]{C20}, the minimax problem \eqref{problem:minimax} can be equivalently casted as the following constrained convex minimization:
\begin{align}\label{problem:minimax2}
	\begin{aligned}
		& \underset{M \in \mathcal{L}(\pi), r}{\min}
		& & r \\
		& \text{s.t.}
		& & D_f(M||L_i) \leq r, ~\textrm{for all}\, i \in \llbracket n \rrbracket.
	\end{aligned}
\end{align}
This problem can be interpreted geometrically as the Chebyshev center problem in the context of Markov chains. For a given $M \in \mathcal{L}(\pi)$ and $r \geq 0$, the set $\{L \in \mathcal{L};~ D_f(M||L) \leq r\}$ can be understood as the ball of Markov generators within radius $r$ of the center $M$. The constraint in \eqref{problem:minimax2} thus entails that this ball encloses the set $\{L_i\}_{i=1}^n$. As the minimization is with respect to the radius $r$, a so-called Chebyshev center is a center of the minimum radius ball that contains $\{L_i\}_{i=1}^n$, and its optimal value is referred as the Chebyshev radius. Precisely we define the Chebyshev center and radius as follows:
\begin{definition}[Chebyshev center and radius]\label{def:Chebyshevcr}
	Consider the minimax problem \eqref{problem:minimax} and its reformulation \eqref{problem:minimax2}. An optimizer 
	$$\argmin_{M \in \mathcal{L}(\pi)} \max_{i \in \llbracket n \rrbracket} D_f(M||L_i)$$
	is called a Chebyshev center with respect to $(f,\{L_i\}_{i=1}^n,\pi)$, and the minimax value $\overline{v}^f$ is called the Chebyshev radius with respect to $(f,\{L_i\}_{i=1}^n,\pi)$.
\end{definition}
For the problem \eqref{problem:minimax2}, we denote the Lagrangian function $\mathsf{L}: \mathbb{R}_+ \times \mathcal{L}(\pi) \times \mathbb{R}_+^n \to \mathbb{R}$ to be
$$\mathsf{L}(r,M,\mathbf{w}) := r + \sum_{i=1}^n w_i(D_f(M||L_i) - r),$$
where $\mathbf{w}$ is the associated Lagrange multiplier. Differentiating $\mathsf{L}$ with respect to $r$ and setting it to zero gives $\mathbf{w}\in \mathcal{S}_n$, the probability simplex that we introduce in \eqref{def:probsim}. The dual problem of \eqref{problem:minimax2} is now written as
\begin{align}
	\max_{\mathbf{w} \in \mathbb{R}_+^n} \min_{r \geq 0, M \in \mathcal{L}(\pi)} \mathsf{L}(r,M,\mathbf{w}) &= 	\max_{\mathbf{w} \in \mathcal{S}_n} \min_{M \in \mathcal{L}(\pi)} \sum_{i=1}^n w_i D_f(M||L_i) \nonumber \\
	&= \max_{\mathbf{w} \in \mathcal{S}_n} \sum_{i=1}^n w_i D_f(M^f_n||L_i), \label{problem:minimax2dual}
\end{align}
where we recall that $M^f_n = M^{f}_n(\mathbf{w},\{L_i\}_{i=1}^n,\pi)$ is the $\mathbf{w}$-weighted information centroid as in Theorem \ref{thm:existuniquecentroid} in the second equality. Note that by weak duality we thus have $\overline{v}^f \geq \max_{\mathbf{w} \in \mathcal{S}_n} \sum_{i=1}^n w_i D_f(M^f_n||L_i)$. Our main result below shows that in fact the strong duality holds and connects various notions introduced earlier on. The proof is stated in Section \ref{subsec:pfmainpure}.

\begin{theorem}\label{thm:mainpure}
	\begin{enumerate}
		\item(Strong duality holds for \eqref{problem:minimax2})\label{it:sd} The strong duality holds for \eqref{problem:minimax2} and there exists $\mathbf{w}^f = \mathbf{w}^f(\{L_i\}_{i=1}^n,\pi) = (w^f_i)_{i=1}^n \in \mathcal{S}_n$ such that the optimal value for the dual problem \eqref{problem:minimax2dual} is attained at $\mathbf{w}^f$, that is,
		$$\overline{v}^f = \sum_{i=1}^n w^f_i D_f(M^{f}_n(\mathbf{w}^f,\{L_i\}_{i=1}^n,\pi)||L_i).$$
		The primal problem \eqref{problem:minimax2} optimal value is attained at $(M^{f}_n(\mathbf{w}^f,\{L_i\}_{i=1}^n,\pi),\overline{v}^f)$, and the optimal value of the problem \eqref{problem:minimax} is attained at $(M^{f}_n(\mathbf{w}^f,\{L_i\}_{i=1}^n,\pi),j)$, where $j \in \llbracket n \rrbracket$ is such that $w^f_j > 0$.
		
		\item(Chebyshev center is a weighted information centroid $M^f_n$)\label{it:ccic} A given pair $(M,r) \in \mathcal{L}(\pi) \times \mathbb{R}_+$ minimizes the primal problem \eqref{problem:minimax2} and a given $\mathbf{w} \in \mathcal{S}_n$ maximizes the dual problem \eqref{problem:minimax2dual} if and only if they satisfy the following two conditions:
		\begin{enumerate}
			\item (Complementary slackness) For $i \in \llbracket n \rrbracket$, we have
			\begin{align*}
				D_f(M||L_i)  \begin{cases}
					 = r  & \text{ if } w_i > 0,\\
					 \leq r & \text{ if } w_i = 0.
				\end{cases}
			\end{align*}
			Furthermore this implies $r = \overline{v}^f$.
			
			\item $M = M^{f}_n(\mathbf{w},\{L_i\}_{i=1}^n,\pi) \in \argmin_{M \in \mathcal{L}(\pi)} \mathsf{L}(r,M,\mathbf{w})$.
		\end{enumerate}
		
		\item(Concavity of the Lagrangian dual)\label{it:cl} The mapping
		$$\mathcal{S}_n \ni \mathbf{w} = (w_i)_{i=1}^n \mapsto \sum_{i=1}^n w_i D_f(M^{f}_n(\mathbf{w},\{L_i\}_{i=1}^n,\pi)||L_i)$$
		is concave.
		
		
		\item(Uniqueness of Chebyshev center under strict convexity of $f$)\label{it:uniquecc}
		Suppose that the parameters $(f,\{L_i\}_{i=1}^n,\pi)$ satisfy the assumptions in Theorem \ref{thm:existuniquecentroid}. If both $\mathbf{w}_1, \mathbf{w}_2 \in \mathcal{S}_n$ maximize the dual problem \eqref{problem:minimax2dual}, we have 
		$$M^{f}_n(\mathbf{w}_1,\{L_i\}_{i=1}^n,\pi) = M^{f}_n(\mathbf{w}_2,\{L_i\}_{i=1}^n,\pi).$$
		In other words, the Chebyshev center is unique.
		
		\item(Characterization of $\overline{v}^f = \underline{v}^f$)\label{it:vfequals} The saddle point property with respect to $(f,\{L_i\}_{i=1}^n,\pi)$ (Definition \ref{def:minimaxvalues}) holds, that is,
		$$\overline{v}^f = \underline{v}^f$$ 
		if and only if there exists $j \in \llbracket n \rrbracket$ and $\mathbf{w}^f = (w^f_i)_{i=1}^n$ a maximizer of the dual problem \eqref{problem:minimax2dual} with $w^f_j > 0$ such that
		\begin{align*}
			M^{f}_n(\mathbf{w}^f,\{L_i\}_{i=1}^n,\pi) \in \argmin_{M \in \mathcal{L}(\pi)} D_f(M||L_j).
		\end{align*}
		In particular, if the parameters $(f,\{L_i\}_{i=1}^n,\pi)$ satisfy the assumptions in Theorem \ref{thm:existuniquecentroid}, then the above statement is equivalent to 
		\begin{align*}
			M^{f}_n(\mathbf{w}^f,\{L_i\}_{i=1}^n,\pi) = M^f(L_j,\pi).
		\end{align*}
		
		\item(Saddle point)\label{it:saddlept} If the saddle point property with respect to $(f,\{L_i\}_{i=1}^n,\pi)$ holds, then 
		$$(M^f(L_l,\pi), L_l)$$
		is a saddle point with respect to $(f,\{L_i\}_{i=1}^n,\pi)$, where $l = \argmax_{i \in \llbracket n \rrbracket} D_f(M^f(L_i,\pi) || L_i)$.
	\end{enumerate}
\end{theorem}

\begin{rk}
	The concavity of the Lagrangian dual in item \eqref{it:cl} plays an important role in developing a projected subgradient algorithm for solving \eqref{problem:minimax2dual} in Section \ref{subsec:algo}.
\end{rk}

In Theorem \ref{thm:mainpure}, we have been investigating the case where $\mathcal{B} = \{L_i\}_{i=1}^n$ in Definition \ref{def:minimaxvalues}, which can be generalized to the convex hull of these generators. Precisely, we shall take
$$\mathcal{B} = \mathrm{conv}((L_i)_{i=1}^n) := \bigg\{\sum_{i=1}^n \alpha_i L_i;~ (\alpha_i)_{i=1}^n \in \mathcal{S}_n, L_i \in \mathcal{L} \text{ for all } i \in \llbracket n \rrbracket \bigg\}$$
and consider minimax problem of the form
\begin{align}\label{problem:minimaxconvexh}
	\inf_{M \in \mathcal{L}(\pi)} \sup_{L \in \mathrm{conv}((L_i)_{i=1}^n)} D_f(M||L) = \overline{v}^f(\mathrm{conv}((L_i)_{i=1}^n),\pi).
\end{align}
An important family of Markov generators that can be written as a convex hull is the family of continuized doubly stochastic Markov generators that we denote by $\mathcal{D}$. We say that $L \in \mathcal{D} \subseteq \mathcal{L}$ if all the row and column sums of $L$ are $0$ and $L = P - I$ where $P$ is a Markov matrix and $I$ is the identity matrix on $\mathcal{X}$, and it can be shown that the discrete uniform distribution on $\mathcal{X}$ is the stationary distribution of such $L$. Let $P_i$ be a permutation matrix on $\mathcal{X}$, then by the Birkhoff-von Neumann theorem \cite[Theorem $8.7.2$]{HJ13} we have
\begin{align}\label{eq:doublysto}
	\mathcal{D} = \mathrm{conv}((P_i-I)_{i=1}^{|\mathcal{X}|!}).
\end{align}
Another family of Markov generators that can be casted under this framework is the set of uniformizable (see \cite[Chapter $2$]{K79}) and $\mu$-reversible generators, which are used in finite truncation of countably infinite Markov chains in queueing theory \citep{V18}, see Example \ref{ex:uniformmix} below.

Our next proposition shows that the analysis in this case of convex hull of generators can be reduced to the setup of Theorem \ref{thm:mainpure}. The proof can be found in Section \ref{subsubsec:pfreducefinite}.
\begin{proposition}[Reduction to the finite case]\label{prop:reducefinite}
	\begin{align*}
		\overline{v}^f(\mathrm{conv}((L_i)_{i=1}^n),\pi) &= \overline{v}^f(\{L_i\}_{i=1}^n,\pi), \\
		\underline{v}^f(\mathrm{conv}((L_i)_{i=1}^n),\pi) &= \underline{v}^f(\{L_i\}_{i=1}^n,\pi).
	\end{align*}
\end{proposition}

The next result collects and combines both Theorem \ref{thm:mainpure} and Proposition \ref{prop:reducefinite}. Its proof is omitted as it is simply restating these two results in a single Corollary.

\begin{corollary}\label{cor:convexhull}
	\begin{enumerate}
	\item(Attainment of the minimax problem \eqref{problem:minimaxconvexh}) There exists $\mathbf{w}^f = \mathbf{w}^f(\{L_i\}_{i=1}^n,\pi) = (w^f_i)_{i=1}^n \in \mathcal{S}_n$ such that the problem \eqref{problem:minimaxconvexh} optimal value is attained at $(M^{f}_n(\mathbf{w}^f,\{L_i\}_{i=1}^n,\pi),L_j)$, where $j \in \llbracket n \rrbracket$ is such that $w^f_j > 0$.
	
	\item(Characterization of $\overline{v}^f = \underline{v}^f$) The saddle point property with respect to $(f,\mathrm{conv}((L_i)_{i=1}^n),\pi)$ (Definition \ref{def:minimaxvalues}) holds, that is,
	$$\overline{v}^f(\mathrm{conv}((L_i)_{i=1}^n),\pi) = \underline{v}^f(\mathrm{conv}((L_i)_{i=1}^n),\pi)$$ 
	if and only if there exists $j \in \llbracket n \rrbracket$ and $\mathbf{w}^f = (w^f_i)_{i=1}^n$ a maximizer of the dual problem \eqref{problem:minimax2dual} with $w^f_j > 0$ such that
	\begin{align*}
		M^{f}_n(\mathbf{w}^f,\{L_i\}_{i=1}^n,\pi) \in \argmin_{M \in \mathcal{L}(\pi)} D_f(M||L_j).
	\end{align*}
	In particular, if the parameters $(f,\{L_i\}_{i=1}^n,\pi)$ satisfy the assumptions in Theorem \ref{thm:existuniquecentroid}, then the above statement is equivalent to 
	\begin{align*}
		M^{f}_n(\mathbf{w}^f,\{L_i\}_{i=1}^n,\pi) = M^f(L_j,\pi).
	\end{align*}
	
	\item(Saddle point) If the saddle point property with respect to $(f,\mathrm{conv}((L_i)_{i=1}^n),\pi)$ holds, then 
	$$(M^f(L_l,\pi), L_l)$$
	is a saddle point with respect to $(f,\mathrm{conv}((L_i)_{i=1}^n),\pi)$, where $l = \argmax_{i \in \llbracket n \rrbracket} D_f(M^f(L_i,\pi),L_i)$.
\end{enumerate}	
	In particular, the above results hold for the set of continuized doubly stochastic Markov generators $\mathcal{D}$ \eqref{eq:doublysto} with $L_i = P_i$ and $n = |\mathcal{X}|!$.
\end{corollary}

\subsection{Results on the pure strategy game}\label{subsec:puregame}

In this subsection, we provide a game-theoretic interpretation of the max-min inequalities and the saddle point property (Definition \ref{def:minimaxvalues}) in the context of the two-person pure strategy game as introduced in Section \ref{subsec:introgame}. 

The max-min inequalities in \eqref{eq:minimaxineq1} and \eqref{eq:minimaxineq2} state that it is advantageous for a player to play second, or more precisely, to know the strategy of the opponent. The difference $\overline{v}^f - \underline{v}^f \geq 0$ can be interpreted as the advantage conferred to a player in knowing the opponent's strategy. If the minimax equality holds, that is, if $\overline{v}^f = \underline{v}^f$, then there is no advantage in playing second or in knowing the strategy of the opponent. 


\textcolor{black}{The main result below characterizes the existence and uniqueness of a Nash equilibrium in this pure strategy game.} The proof is deferred to Section \ref{subsubsec:pfpureNash}.

\begin{corollary}\label{cor:pureNash}
	\textcolor{black}{Consider the two-person pure strategy game as discussed in Section \ref{subsec:introgame}.} Let $\mathcal{B}$ be either $\{L_i\}_{i=1}^n$ or its convex hull $\mathrm{conv}((L_i)_{i=1}^n)$.
	\begin{enumerate}
		\item(Characterization of pure strategy Nash equilibrium)\label{it:existNash} A pure strategy Nash equilibrium with respect to $(f,\mathcal{B},\pi)$  exists
		if and only if there exists $j \in \llbracket n \rrbracket$ and $\mathbf{w}^f = (w^f_i)_{i=1}^n$ a maximizer of the dual problem \eqref{problem:minimax2dual} with $w^f_j > 0$ such that
		\begin{align*}
			M^{f}_n(\mathbf{w}^f,\{L_i\}_{i=1}^n,\pi) \in \argmin_{M \in \mathcal{L}(\pi)} D_f(M||L_j).
		\end{align*}
		In particular, if the parameters $(f,\{L_i\}_{i=1}^n,\pi)$ satisfy the assumptions in Theorem \ref{thm:existuniquecentroid}, then the above statement is equivalent to 
		\begin{align*}
			M^{f}_n(\mathbf{w}^f,\{L_i\}_{i=1}^n,\pi) = M^f(L_j,\pi).
		\end{align*}
		
		\item(Uniqueness of pure strategy Nash equilibrium)\label{it:uniqueNash} Suppose that the parameters $(f,\{L_i\}_{i=1}^n,\pi)$ satisfy the assumptions in Theorem \ref{thm:existuniquecentroid} and a pure strategy Nash equilibrium exists, which is given by 
		$$(M^f(L_l,\pi), L_l),$$
		where $l \in \argmax_{i \in \llbracket n \rrbracket} D_f(M^f(L_i,\pi),L_i)$.
		It is unique if and only if the index $$l = \argmax_{i \in \llbracket n \rrbracket} D_f(M^f(L_i,\pi),L_i)$$ is unique.
	\end{enumerate}	
\end{corollary}

\subsubsection{Examples}\label{subsubex:pureex}

In this subsection, we give a few simple yet illustrative examples to demonstrate the theory that we have developed thus far. We shall see that, depending on the parameters of the game $(f,\{L_i\}_{i=1}^n,\pi)$, the associated pure strategy Nash equilibrium or saddle point may or may not exist.

\begin{example}[An example with multiple pure strategy Nash equilibria (or saddle points)]\label{ex:1pure}
	In the first example, we take $n = 2$ and consider two generators with $L_1 = L$ and $L_2 = L_{\pi}$, where $L \in \mathcal{L}$ and recall that $L_{\pi}$ is the $\pi$-dual of $L$.
	
	In view of the bisection property \citep{CW23}, we have $D_f(M||L_1) = D_f(M||L_2)$, and hence for any weight $\mathbf{w} \in \mathcal{S}_2$, we see that
	$$M^{f}_2(\mathbf{w},\{L_i\}_{i=1}^2,\pi) \in \argmin_{M \in \mathcal{L}(\pi)} D_f(M||L_i)$$
	for $i = 1,2$. Thus, a pure strategy Nash equilibrium or saddle point, with respect to $(f,\{L_i\}_{i=1}^2,\pi)$, exists, according to Corollary \ref{cor:pureNash}.
	
	Now, we further assume that the parameters $(f,\{L_i\}_{i=1}^n,\pi)$ satisfy the assumptions in Theorem \ref{thm:existuniquecentroid}. Using the bisection property again, the two pure strategy Nash equilibria are
	$$(M^f(L,\pi),L), \quad (M^f(L,\pi),L_{\pi}).$$
\end{example}

\begin{example}[An example with a unique pure strategy Nash equilibrium (or saddle point)]
	In the second example, we again take $n = 2$ and consider two generators with $L_1, L_2 \in \mathcal{L}$. We further assume that $D_f$ is the $\alpha$-divergence, and $L_1,L_2$ are chosen to satisfy the condition that
	\begin{align}\label{eq:examplecond}
		D_f(M^f(L_1,\pi)||L_1) > D_f(M^f(L_1,\pi)||L_2).
	\end{align}
	Applying the Pythagorean identity of the $\alpha$-divergence \citep{CW23} to the right hand side above, we thus have
	 \begin{align}\label{eq:examplecond2}
	 	D_f(M^f(L_1,\pi)||L_1) > D_f(M^f(L_1,\pi)||L_2) &= D_f(M^f(L_2,\pi)||L_2) +  D_f(M^f(L_1,\pi)||M^f(L_2,\pi)) \nonumber \\
	 	&\geq D_f(M^f(L_2,\pi)||L_2).
	 \end{align}
 
 	Now, we claim that in this example, the pair $(M^f(L_1,\pi),D_f(M^f(L_1,\pi)||L_1)) \in \mathcal{L}(\pi) \times \mathbb{R}_+$ minimizes the primal problem \eqref{problem:minimax2} and $\mathbf{w} = (1,0) \in \mathcal{S}_2$ maximizes the dual problem \eqref{problem:minimax2dual}. This readily follows from item \eqref{it:ccic} in Theorem \ref{thm:mainpure}, where the condition \eqref{eq:examplecond} ensures that the complementary slackness holds.
 	
 	Using item \eqref{it:existNash} in Corollary \ref{cor:pureNash}, we see that a Nash equilibrium exists, and in view of \eqref{eq:examplecond2} and item \eqref{it:uniqueNash} in Corollary \ref{cor:pureNash}, the unique pure strategy Nash equilibrium or saddle point is given by
 	$$(M^f(L_1,\pi),L_1).$$
\end{example}

\begin{example}[An example where a pure strategy Nash equilibrium (or saddle point) fails to exist]\label{ex:3pure}
	In the final example, we again specialize into $n = 2$ and take two Markov generators with $L_1, L_2 \in \mathcal{L}$. We further assume that $D_f$ is the $\alpha$-divergence, and $L_1,L_2$ are chosen to satisfy the conditions that $M^f(L_1,\pi) \neq M^f(L_2,\pi)$ and
	\begin{align}\label{eq:examplecond3}
		D_f(M^f(L_1,\pi)||L_1) = D_f(M^f(L_2,\pi)||L_2).
	\end{align}
	
	Suppose \textcolor{black}{for contradiction} that a pure strategy Nash equilibrium exists. In view of Corollary \ref{cor:pureNash} and $f$ is strictly convex under $\alpha$-divergence, the only possible candidates for the maximizer of the dual problem is either $\mathbf{w}^f = (1,0)$ or $\mathbf{w}^f = (0,1)$.
	
	In the former case, we thus have $D_f(M^f(L_1,\pi)||L_1) = r$ and
	\begin{align*}
		r &= D_f(M^f(L_1,\pi)||L_1) \\
		  &\geq D_f(M^f(L_1,\pi)||L_2) \\
		  &=  D_f(M^f(L_2,\pi)||L_2) +  D_f(M^f(L_1,\pi)||M^f(L_2,\pi)) \\
		  &=  D_f(M^f(L_1,\pi)||L_1) +  D_f(M^f(L_1,\pi)||M^f(L_2,\pi)),
	\end{align*}
	where we use the Pythagorean identity \citep{CW23} in the second equality and \eqref{eq:examplecond3} in the last equality. This implies that $D_f(M^f(L_1,\pi)||M^f(L_2,\pi)) = 0$ and hence $M^f(L_1,\pi) = M^f(L_2,\pi)$, which contradicts the assumption. Analogously we can handle the case of $\mathbf{w}^f = (0,1)$.
	
	As a result, we see that a pure strategy Nash equilibrium fails to exist. For further discussions on this example, please refer to Example \ref{ex:3mix}.
\end{example}

\subsection{Results on the mixed strategy game}\label{subsec:mixgame}

\textcolor{black}{This subsection concerns the two-person mixed strategy game as introduced in Section \ref{subsec:introgame}.} Before we state the main results of this subsection, we give an important proposition that connects various minimax and maximin values we have introduced thus far, namely $\overline{v}^f,\overline{V}^f, \underline{V}^f$ and the dual problem \eqref{problem:minimax2dual}. The proof is deferred to Section \ref{subsubsec:pfminimaxineq}.

\begin{proposition}\label{prop:minimaxeqmixed}
	For any $\mathcal{B}\subseteq \mathcal{L}$, we have
	\begin{align}\label{eq:minimaxeqmixed1}
		\overline{v}^f(\mathcal{B},\pi) \geq \overline{V}^f(\mathcal{B},\pi) \geq \underline{V}^f(\mathcal{B},\pi),
	\end{align}
	where we recall that $\overline{v}^f$ is introduced in Definition \ref{def:minimaxvalues} while $\overline{V}^f,\underline{V}^f$ are defined in Definition \ref{def:minimaxmix}. In particular, if $\mathcal{B}$ is either $\{L_i\}_{i=1}^n$ or its convex hull $\mathrm{conv}((L_i)_{i=1}^n)$, this leads to
	\begin{align}\label{eq:minimaxeqmixed2}
		\overline{v}^f(\mathcal{B},\pi) \geq \overline{V}^f(\mathcal{B},\pi) \geq \underline{V}^f(\mathcal{B},\pi) \geq \max_{\mathbf{w} \in \mathcal{S}_n} \sum_{i=1}^n w_i D_f(M^f_n||L_i).
	\end{align}
\end{proposition}

In the setting of \eqref{eq:minimaxeqmixed2}, by Theorem \ref{thm:mainpure} item \eqref{it:sd}, the strong duality holds and hence we have
$$\overline{v}^f(\mathcal{B},\pi) = \max_{\mathbf{w} \in \mathcal{S}_n} \sum_{i=1}^n w_i D_f(M^f_n||L_i).$$
This forces a minimax equality, that is,
$$\overline{V}^f(\mathcal{B},\pi) = \underline{V}^f(\mathcal{B},\pi).$$
In other words, all inequalities are in fact equalities in \eqref{eq:minimaxeqmixed2}. We collect this result and a few others in Corollary \ref{cor:convexhull} into the following main results:

\begin{theorem}\label{thm:mainmix}
	Consider the two-person mixed strategy game with respect to the parameters $(f, \mathcal{B},\pi)$, where $\mathcal{B}$ is either $\{L_i\}_{i=1}^n$ or its convex hull $\mathrm{conv}((L_i)_{i=1}^n)$. Recall that there exists $\mathbf{w}^f = \mathbf{w}^f(\{L_i\}_{i=1}^n,\pi) = (w^f_i)_{i=1}^n \in \mathcal{S}_n$ such that the optimal value for the dual problem \eqref{problem:minimax2dual} is attained at $\mathbf{w}^f$ by Theorem \ref{thm:mainpure}.
	\begin{enumerate}
		\item(A mixed strategy Nash equilibrium always exists) The minimax equality holds, that is,
		$$\overline{V}^f(\mathcal{B},\pi) = \underline{V}^f(\mathcal{B},\pi) = \sum_{i=1}^n w^f_i D_f(M^{f}_n(\mathbf{w}^f,\{L_i\}_{i=1}^n,\pi)||L_i) =: V,$$
		and hence $V$ is the value of the mixed strategy game. The pair $(M^{f}_n(\mathbf{w}^f,\{L_i\}_{i=1}^n,\pi),\mathbf{w}^f) \in \mathcal{L}(\pi) \times \mathcal{P}(\mathcal{B})$ is a mixed strategy Nash equilibrium with respect to the parameters $(f, \mathcal{B},\pi)$.
		
		\item(Game-theoretic interpretation of Chebyshev center and Chebyshev radius) Recall the definition of Chebyshev center and Chebyshev radius in Definition \ref{def:Chebyshevcr}. A minimax strategy of Nature is given by $M^{f}_n(\mathbf{w}^f,\{L_i\}_{i=1}^n,\pi)$, a weighted information centroid, which is also a Chebyshev center. 
		In view of the complementary slackness in Theorem \ref{thm:mainpure}, we have
		$$V = D_f(M^{f}_n(\mathbf{w}^f,\{L_i\}_{i=1}^n,\pi)||L_l),$$
		where the index $l$ satisfies $w^f_l > 0$. In words, the Chebyshev radius is the value of the game $V$.
				
		If the parameters $(f,\{L_i\}_{i=1}^n,\pi)$ satisfy the assumptions in Theorem \ref{thm:existuniquecentroid}, then $M^{f}_n(\mathbf{w}^f,\{L_i\}_{i=1}^n,\pi)$ is the unique minimax strategy.
		
		\item(Bayes risk with respect to $\mathbf{w}^f$ is value of the game) $\mathbf{w}^f$ is a maximin strategy of the probabilist. $M^{f}_n(\mathbf{w}^f,\{L_i\}_{i=1}^n,\pi)$ is a Bayes strategy with respect to $\mathbf{w}^f$, and the Bayes risk with respect to $\mathbf{w}^f$ is the value of the game $V$.
	\end{enumerate}
\end{theorem}

\begin{rk}
	It is interesting to note that a mixed strategy Nash equilibrium always exists in the two-person mixed strategy game (Theorem \ref{thm:mainmix}) while a pure strategy Nash equilibrium may or may not exist in the two-person pure strategy game (Corollary \ref{cor:pureNash}). For an example of the latter case we recall Example \ref{ex:3pure} and another Example \ref{ex:3mix} below. This is analogous to the classical two-person zero-sum game in the game theory literature where a pure strategy Nash equilibrium may not exist, see \cite[Chapter $2$]{KarlinPeres17}.
\end{rk}

\begin{rk}[Game-theoretic consequences of the mixed strategy game]
	As a mixed strategy Nash equilibrium always exists, there is no advantage for a player (the probabilist or Nature) to play second, or more generally, to know the strategy of the opponent in the mixed strategy game. 
\end{rk}

\begin{rk}[Another proof of $\overline{V}^f(\mathcal{B},\pi) = \underline{V}^f(\mathcal{B},\pi)$ via the Sion's minimax theorem]\label{rk:Sion}
	One common strategy in establishing minimax results, in the context of source coding and information theory, relies on the Sion's minimax theorem, see for instance \cite[Theorem $35$]{EH14}. For given $\mu \in \mathcal{P}(\mathcal{B})$ and $M \in \mathcal{L}(\pi)$, the mapping
	$$(M,\mu) \mapsto \int_{\mathcal{B}} D_f(M||L)\, \mu(dL)$$
	is clearly concave in $\mu$ and convex in $M$. As we are minimizing over $M \in \mathcal{L}(\pi)$, the set of $\pi$-reversible generator $\mathcal{L}(\pi)$ is convex yet unbounded, and hence it is not a compact subset of $\mathcal{L}$. On the other hand, as $\mathcal{B}$ is either a finite set of Markov generators $\{L_i\}_{i=1}^n$ or its convex hull $\mathrm{conv}((L_i)_{i=1}^n)$, the set $\mathcal{P}(\mathcal{B})$ is thus compact. As a result, the Sion's minimax theorem is readily applicable in this setting which yields
	$$\overline{V}^f(\mathcal{B},\pi) = \underline{V}^f(\mathcal{B},\pi).$$
	The proof that we presented in \eqref{eq:minimaxeqmixed2} relies on the Lagrangian duality theory as in Theorem \ref{thm:mainpure}, which naturally gives fine properties concerning the mixed strategy Nash equilibrium and related important objects such as $\mathbf{w}^f$ and connects with earlier sections of the manuscript. These properties do not seem to be immediate consequences or corollaries from the Sion's minimax theorem.
\end{rk}

\subsubsection{Examples}\label{subsubex:mixedex}

In this subsection, similar to Section \ref{subsubex:pureex}, we detail four simple examples to demonstrate the theory concerning the two-person mixed strategy game.

\begin{example}[Answering a question of Laurent Miclo: an example with the existence of multiple mixed strategy Nash equilibria]\label{ex:1mix}
	In the first example, we continue the setting discussed in Example \ref{ex:1pure}. That is, we consider two ($n = 2$) generators with $L_1 = L$ and $L_2 = L_{\pi}$, where $L \in \mathcal{L}$.
	
    Using the bisection property of $D_f$ \cite{CW23}, we have $D_f(M||L_1) = D_f(M||L_2)$, and for any weight $\mathbf{w} \in \mathcal{S}_2$, we see that
	$$M^{f}_2(\mathbf{w},\{L_i\}_{i=1}^2,\pi) \in \argmin_{M \in \mathcal{L}(\pi)} D_f(M||L_i)$$
	for $i = 1,2$. In other words, any pair of the form $(M^f(L_1,\pi),\mathbf{w})$ is a mixed strategy Nash equilibrium with respect to the parameters of this game. Note that this game also has multiple pure strategy Nash equilibria as shown in Example \ref{ex:1pure}.
	
	In the special case of considering $f(t) = |t-1|$, the $f$-divergence is the total variation distance. We also recall $P_{-\infty}, P_{\infty}$ are introduced in \eqref{eq:P-infty} and \eqref{eq:Pinfty} respectively. It is shown in \cite{CH18} that any convex combinations of $P_{-\infty}$ and $P_{\infty}$ minimize
	$$\{a P_{-\infty} + (1-a) P_{\infty};~ a \in [0,1]\} \subseteq \argmin_{M \in \mathcal{L}(\pi)} D_f(M||L_i).$$
	As such, based on our earlier analysis, $(P_{\infty},\mathbf{w})$ and $(P_{-\infty},\mathbf{w})$ are mixed strategy Nash equilibria of the game under the total variation distance. This gives a possible answer to a question of Laurent Miclo on game-theoretic interpretations of $P_{-\infty}, P_{\infty}$.
	
	If $f$ is strictly convex with $f^{\prime}(1) = 0$, for instance the $\alpha$-divergence, then $M^f(L_1,\pi)$ is unique, and hence it is the unique minimax strategy of Nature in the mixed strategy game, while any $\mathbf{w} \in \mathcal{S}_2$ is a maximin strategy of the probabilist. $M^f(L_1,\pi)$ is also the unique Bayes strategy with respect to any $\mathbf{w}$. The value of the game is given by $V = D_f(M^f(L_1,\pi)||L_1)$. 
	
	Note that $M^f(L_1,\pi)$ is also the unique minimax strategy of Nature in the pure strategy game.
\end{example}

%
%

\begin{example}[An example where a pure strategy Nash equilibrium fails to exist and a mixed strategy Nash equilibrium (always) exists]\label{ex:3mix}
	In the second example, we continue our investigation in Example \ref{ex:3pure} and take $n = 2$ along with two Markov generators where $L_1, L_2 \in \mathcal{L}$. We choose $D_f$ to be the $\alpha$-divergence, and $L_1,L_2$ are assumed to satisfy the conditions that $M^f(L_1,\pi) \neq M^f(L_2,\pi)$ and
	\begin{align}\label{eq:examplecond4}
		D_f(M^f(L_1,\pi)||L_1) = D_f(M^f(L_2,\pi)||L_2).
	\end{align}
	
	Recall that in Example \ref{ex:3pure}, we have already shown that a pure strategy Nash equilibrium cannot exist with these choices of parameters, and furthermore a maximizer of the dual problem \eqref{problem:minimax2dual}, which exists, cannot be $\mathbf{w}^f = (1,0)$ and $\mathbf{w}^f = (0,1)$.
	
	We thus come to the conclusion that $w^f_i > 0$ for $i = 1,2$. The weighted information centroid $M^{f}_2(\mathbf{w}^f,\{L_i\}_{i=1}^2,\pi)$ is the unique minimax strategy of Nature, and is also the unique Bayes strategy with respect to $\mathbf{w}^f$.
	
	Figure \ref{fig:purenotexist} demonstrates a possible visualization of this example.
	
	\begin{figure}[h]
		\centering
		\includegraphics[width=0.8\textwidth]{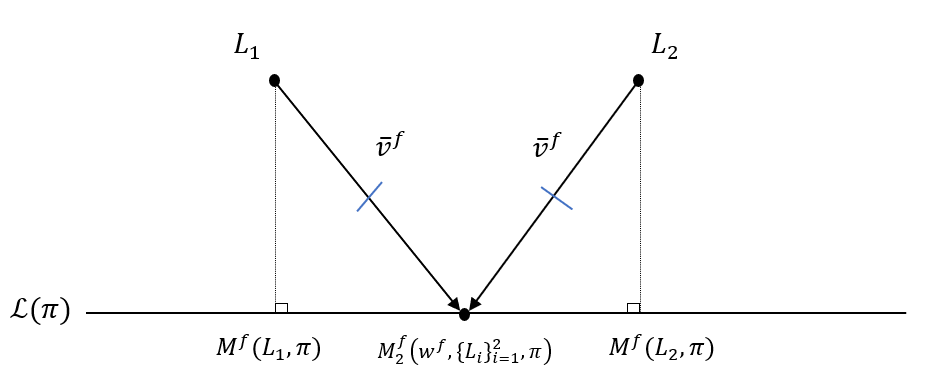}
		\caption{A visualization of Example \ref{ex:3mix}. Note that the mixed strategy Nash equilibrium is $(M^{f}_2(\mathbf{w}^f,\{L_i\}_{i=1}^2,\pi),\mathbf{w}^f)$, while a pure strategy Nash equilibrium does not exist. The complementary slackness in Theorem \ref{thm:mainpure} ensures that $D_f(M^{f}_2(\mathbf{w}^f,\{L_i\}_{i=1}^2,\pi)||L_1) = D_f(M^{f}_2(\mathbf{w}^f,\{L_i\}_{i=1}^2,\pi)||L_2) = \overline{v}^f = V$, the value of the mixed strategy game. \eqref{eq:examplecond4} also holds in this figure.}
		\label{fig:purenotexist}
	\end{figure}
\end{example}

\begin{example}[An example where a pure strategy Nash equilibrium and mixed strategy Nash equilibrium coincides]\label{ex:puremixsame}
	In this example, we take $n=3$ Markov generators $L_i \in \mathcal{L}$ for $i \in \llbracket 3 \rrbracket$, and they are chosen in a special way such that
	$$D_f(M^f(L_2,\pi)||L_1) = D_f(M^f(L_2,\pi)||L_3) = D_f(M^f(L_2,\pi)||L_2).$$
	At the same time, we require that, for $j \in \{1,3\}$,
	$$D_f(M^f(L_2,\pi)||L_2) > D_f(M^f(L_j,\pi)||L_j).$$
	
	With these parameter choices of the game, using Theorem \ref{thm:mainpure} and \ref{thm:mainmix} we see that we can take $\mathbf{w}^f = (0,1,0)$, and $(M^f(L_2,\pi),\mathbf{w}^f)$ is a mixed strategy Nash equilibrium while {\color{black}$(M^f(L_2,\pi),L_2)$} is a pure strategy Nash equilibrium. Figure \ref{fig:puremixsame} provides a visualization in this setting.
	
	\begin{figure}[h]
		\centering
		\includegraphics[width=0.8\textwidth]{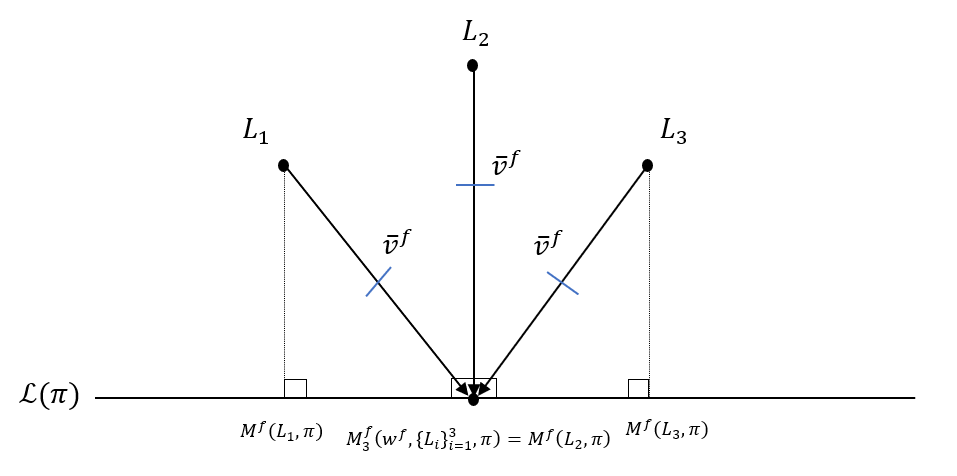}
		\caption{A visualization of Example \ref{ex:puremixsame}. }
		\label{fig:puremixsame}
	\end{figure}
\end{example}

\begin{example}[An example of uniformizable and $\mu$-reversible generators]\label{ex:uniformmix}
	Recall that we have been considering settings where $\mathcal{B}$ is either $\{L_i\}_{i=1}^n$ or its convex hull $\mathrm{conv}((L_i)_{i=1}^n)$. This covers for instance $\mathcal{D}$, the set of doubly stochastic Markov generators, as discussed in \eqref{eq:doublysto}. The aim of this example is to show that this also includes the set of uniformizable and $\mu$-reversible generators, where $\mu \neq \pi$.
	
	Let $L \in \mathcal{L}(\mu)$. $L$ is said to be $\lambda$-uniformizable if $\max_{x \in \mathcal{X}} |L(x,x)| \leq \lambda$, and we denote by $\mathcal{L}_{\lambda}(\mu) \subset \mathcal{L}(\mu)$ to be the set of such Markov generators. Without loss of generality, in this example we let $m = |\mathcal{X}|$ and label the state space $\mathcal{X} = \{1,2,\ldots,m\}$. For a given $\lambda > 0$, we now claim that 
	$$\mathcal{L}_{\lambda}(\mu) \subseteq \mathrm{conv}(\{L_{x,y}\}_{x<y,x,y \in \mathcal{X}} \cup \{\mathbf{0}\}),$$
	where for $x < y \in \mathcal{X}$, $L_{x,y}(x,y) := \lambda \dfrac{m(m-1)}{2}$ while $L_{x,y}(y,x) := \lambda \dfrac{m(m-1)}{2} \dfrac{\mu(x)}{\mu(y)}$ and is zero otherwise for other off-diagonal entries of $L_{x,y}$. We also recall that $\mathbf{0}$ is the all-zero Markov generator. Let $L \in \mathcal{L}_{\lambda}(\mu)$. Making use of the $\mu$-reversibility, we write that
	$$L = \sum_{x < y} w_{x,y} L_{x,y} + \left(1-\sum_{x < y} w_{x,y}\right) \mathbf{0},$$
	where $w_{x,y} := \dfrac{L(x,y)}{\lambda \dfrac{m(m-1)}{2}}$.
	
	As a consequence, various results that have been stated are readily applicable to the set $\mathcal{L}_{\lambda}(\mu)$.
	
	
\end{example}

\subsection{A projected subgradient algorithm to find an approximate mixed strategy Nash equilibrium (or Chebyshev center)}\label{subsec:algo}

The aim of this \textcolor{black}{subsection} is to develop a simple and easy-to-implement projected subgradient method to find an approximate mixed strategy Nash equilibrium. Throughout this \textcolor{black}{subsection}, we consider the two-person mixed strategy game introduced in Section \ref{subsec:introgame} with respect to the parameters $(f, \mathcal{B},\pi)$, where $\mathcal{B}$ is either $\{L_i\}_{i=1}^n$ or its convex hull $\mathrm{conv}((L_i)_{i=1}^n)$. Recall that $\mathbf{w}^f$ is a maximin strategy of the probabilist, and the pair $(M^{f}_n(\mathbf{w}^f,\{L_i\}_{i=1}^n,\pi),\mathbf{w}^f)$ is a mixed strategy Nash equilibrium. To find such an equilibrium algorithmically, one crucial step amounts to solving the corresponding dual problem \eqref{problem:minimax2dual}. If the weighted information centroid admits a closed-form expression, for instance in the case of $\alpha$-divergence or the examples given in Theorem \ref{thm:examplecentroid}, then making use of $\mathbf{w}^f$ a mixed strategy Nash equilibrium is found. 

Instead of solving the maximization problem in \eqref{problem:minimax2dual}, we reformulate \textcolor{black}{it} into the following equivalent minimization problem by multiplying by $-1$, in order for us to consider the subgradient (rather than supergradient) method in the optimization literature: 
\begin{align*}
	\min_{\mathbf{w} \in \mathcal{S}_n} h(\mathbf{w})
\end{align*}
where
\begin{align}\label{def:h}
	h(\mathbf{w}) := - \sum_{i=1}^n w_i D_f(M^{f}_n(\mathbf{w},\{L_j\}_{j=1}^n,\pi)||L_i)
\end{align}
Note that the dependence of $h$ on $(f, \mathcal{B},\pi)$ is suppressed to avoid notational burden in the definition above. In view of Theorem \ref{thm:mainpure} item \eqref{it:cl}, the function $h$ is convex.

The algorithm that we propose is novel from two perspectives. First, no algorithm has been developed in the context of \textcolor{black}{the game of the statistician against Nature} \citep{H97,HO97,GZ06} to compute the Nash equilibrium, while in this paper we propose a simple subgradient algorithm to approximately find the mixed strategy equilibrium. Second, our algorithm harnesses on the centroid structure and the subgradient of $h$, which has not been observed in the literature of Chebyshev center computation of probability measures \citep{EBT08,C20}.

In the first main result of this Section, we identify a subgradient of $h$, and its proof is deferred to Section \ref{subsubsec:pfsubgh}:
\begin{theorem}[Subgradient of $h$ and an upper bound of its $\ell^2$-norm]\label{thm:subgh}
	A subgradient of $h$ at $\mathbf{v} \in \mathcal{S}_n$ is given by $\mathbf{g} = \mathbf{g}(\mathbf{v}) = (g_1,g_2,\ldots,g_n) \in \mathbb{R}^n$, where for $i \in \llbracket n \rrbracket$, we have
	\begin{align}\label{eq:subgi}
		g_i = D_f(M^f_n(\mathbf{v},\{L_j\}_{j=1}^n,\pi)||L_n) - D_f(M^f_n(\mathbf{v},\{L_j\}_{j=1}^n,\pi)||L_i).
	\end{align}
	That is, $\mathbf{g}$ satisfies, for $\mathbf{w},\mathbf{v} \in \mathcal{S}_n$,
	$$h(\mathbf{w}) \geq h(\mathbf{v}) + \sum_{i=1}^n g_i (w_i - v_i).$$
	Moreover, the $\ell^2$-norm of $\mathbf{g}(\mathbf{v})$ is bounded above by
	\begin{align}\label{eq:subgnormbd}
		\norm{\mathbf{g}(\mathbf{v})}_2^2 := \sum_{i=1}^n g_i^2 \leq n \left(|\mathcal{X}| \sup_{\mathbf{v} \in \mathcal{S}_n;~i \in \llbracket n \rrbracket;~L_i(x,y)>0} L_i(x,y) f\left(\dfrac{M^f_n(\mathbf{v},\{L_j\}_{j=1}^n,\pi)(x,y)}{L_i(x,y)}\right)\right)^2 =: B. 
	\end{align}
\end{theorem}

\begin{rk}
	The choice of $L_n$ in the first term of $g_i$ in \eqref{eq:subgi} is in fact completely arbitrary, and we shall see in the proof that it can be replaced by any $L_l$ with $l \in \llbracket n \rrbracket$. We thus obtain at least $n$ subgradients of $h$ at a given point $\mathbf{v} \in \mathcal{S}_n$.
\end{rk}

\begin{rk}
	The upper bound $B$ plays an important role in determining the convergence rate of the projected subgradient algorithm, see Theorem \ref{thm:subgconv} and its proof below.
\end{rk}

We proceed to develop and analyze a projected subgradient algorithm to solve the dual problem \eqref{problem:minimax2dual}. Suppose that the number of iterations to be ran is $t$.

For $i = 1,2,\ldots,t$, we first update the weights of the current iteration by its subgradient, that is, we set
\begin{align*}
	\mathbf{v}^{(i)} = \mathbf{w}^{(i-1)} - \eta \mathbf{g}(\mathbf{w}^{(i-1)}),
\end{align*}
where $\eta > 0$ is the constant stepsize of the algorithm while we recall that $\mathbf{g}$ is a subgradient of $h$ as in \eqref{eq:subgi}. In the second step, we project $\mathbf{v}^{(i)}$ onto the probability simplex $\mathcal{S}_n$ by computing
\begin{align*}
	\mathbf{w}^{(i)} = \argmin_{\mathbf{w} \in \mathcal{S}_n} \norm{\mathbf{w} - \mathbf{v}^{(i)}}_2^2.
\end{align*}

This can be done by \textcolor{black}{applying} existing efficient projection algorithms (see e.g. \cite{Condat16}), and we do not further investigate these projection algorithms in this manuscript. We then repeat the above two steps for $i = 1,2,\ldots,t$. The output of the algorithm is the sequence $(\mathbf{w}^{(i)})_{i=1}^t$. Precisely, the algorithm is stated in Algorithm \ref{algo:subg}.

\begin{algorithm}
	\caption{A projected subgradient algorithm to find an approximate mixed strategy Nash equilibrium}\label{algo:subg}
	
	\SetKwInOut{Input}{Input}
	\SetKwInOut{Output}{Output}
	
	\Input{
		Initial weight value $\mathbf{w}^{(0)} \in \mathcal{S}_n$, function $f$ (and $D_f$), set $\{L_i\}_{i=1}^n$, target distribution $\pi$, stepsize $\eta > 0$ and number of iterations $t$.
	}
	\Output{
		The sequence $(\mathbf{w}^{(i)})_{i=1}^t$.
	}
	
	\For{$i = 1,2,\ldots,t$}{
		\tcp{Update via subgradient}
		$\mathbf{v}^{(i)} = \mathbf{w}^{(i-1)} - \eta \mathbf{g}(\mathbf{w}^{(i-1)})$\;
		\tcp{Projection onto $\mathcal{S}_n$}
		$\mathbf{w}^{(i)} = \argmin_{\mathbf{w} \in \mathcal{S}_n} \norm{\mathbf{w} - \mathbf{v}^{(i)}}_2^2.$
	}
\end{algorithm}

Before we proceed to analyze the convergence of Algorithm \ref{algo:subg}, \textcolor{black}{we first build some intuition about why it might possibly work}. Suppose that we are in the hypothetical setting where the optimal weights are all positive, that is, $w^f_i > 0$ for all $i \in \llbracket n \rrbracket$. Recall that by the complementary slackness condition in Theorem \ref{thm:mainpure}, in this setting the subgradient of $h$ at $\mathbf{w}^f$ is zero since
$$g_i(\mathbf{w}^f) = D_f(M^f_n(\mathbf{w}^f,\{L_j\}_{j=1}^n,\pi)||L_n) - D_f(M^f_n(\mathbf{w}^f,\{L_j\}_{j=1}^n,\pi)||L_i) = 0.$$
As a result, once Algorithm \ref{algo:subg} reaches $\mathbf{w}^f$, it stays put and no longer moves again, which is a desirable behaviour as an optimal value has been reached. On the other hand, if at the current iteration the subgradient is $g_i > 0$ at some $i$, then the weights are updated along the subgradient direction followed by projection onto $\mathcal{S}_n$.

The success of Algorithm \ref{algo:subg} relies crucially on computation of the subgradient $\mathbf{g}(\mathbf{w})$ of $h$ at $\mathbf{w}$, which further depends on the expression 
$M^f_n(\mathbf{w},\{L_j\}_{j=1}^n,\pi)$. For example, it is known in the case of $\alpha$-divergence or the examples given in Theorem \ref{thm:examplecentroid}, and hence the algorithm can be readily applied in these settings. Note that in the implementation of Algorithm \ref{algo:subg}, we assume that the centroid $M^f_n(\mathbf{w},\{L_j\}_{j=1}^n,\pi)$ is accessible and we do not discuss algorithms to compute these centroids numerically.

In this paper, we do not consider \textcolor{black}{the many possible} variants of the projected subgradient algorithm in Algorithm \ref{algo:subg}, for instance algorithms with changing or adaptive stepsize, or algorithms with stochastic subgradients.

Our second main result of this \textcolor{black}{subsection} gives a convergence rate of $\mathcal{O}(1/\sqrt{t})$ of Algorithm \ref{algo:subg} with an appropriate choice of stepsize. The proof can be found in Section \ref{subsubsec:pfsubgconv}. We also write
\begin{align*}
	\overline{\mathbf{w}}^t := \dfrac{1}{t}\sum_{i=1}^t \mathbf{w}^{(i)},
\end{align*}
the arithmetic average up to iteration $t$ of the outputs of Algorithm \ref{algo:subg}.

\begin{theorem}\label{thm:subgconv}
	Consider Algorithm \ref{algo:subg} with its output $(\mathbf{w}^{(i)})_{i=1}^t$ and the notations therein. We have
	\begin{align}\label{eq:subgconvbd}
		h\left(\overline{\mathbf{w}}^t\right) - h\left(\mathbf{w}^{f}\right) \leq \dfrac{n}{2\eta t} + \dfrac{\eta B}{2},
	\end{align}
	where we recall that $h$ is the function defined in \eqref{def:h} and $B$ is introduced in \eqref{eq:subgnormbd}. In particular, if we take the constant stepsize to be
	$$\eta = \sqrt{\dfrac{n}{tB}},$$
	we thus have
	$$h\left(\overline{\mathbf{w}}^t\right) - h\left(\mathbf{w}^{f}\right) \leq \dfrac{1}{2}\sqrt{\dfrac{nB}{t}} + \dfrac{1}{2}\sqrt{\dfrac{nB}{t}} = \sqrt{\dfrac{nB}{t}}.$$
	Given an arbitrary $\varepsilon > 0$, if we further choose
	$$t = \left\lceil \dfrac{n B}{\varepsilon^2} \right\rceil,$$
	then we reach an $\varepsilon$-close value to $ h\left(\mathbf{w}^{f}\right)$ in the sense that
	$$h\left(\overline{\mathbf{w}}^t\right) - h\left(\mathbf{w}^{f}\right) \leq \varepsilon.$$
\end{theorem}

Before we introduce the final main result of this Section, we shall fix a few notations. 
We shall consider the set of continuized generators that does not stay at the same state in one step, that is, for $i \in \llbracket n \rrbracket$ we consider
\begin{align}\label{eq:Liclass}
	L_i := P_i - I,
\end{align}
where $I$ is the identity matrix on $\mathcal{X}$ and $P_i$ is a transition matrix with $P_i(x,x) = 0$ for all $x \in \mathcal{X}$. For two Markov generators $L,M \in \mathcal{L}$, by choosing $f(t) = |t-1|/2$, the total variation distance between $L,M$ is given by
$$D_{\mathrm{TV}}(L||M) := D_f(L||M) = \dfrac{1}{2} \sum_{x \in \mathcal{X}} \pi(x) \sum_{y \in \mathcal{X}\backslash\{x\} }|L(x,y) - M(x,y)| = D_{\mathrm{TV}}(M||L).$$

For the family of $\{L_i\}_{i=1}^n$ satisfying \eqref{eq:Liclass}, we quantify the convergence of $M^f_n(\overline{\mathbf{w}}^t,\{L_j\}_{j=1}^n,\pi)$ towards $M^f_n(\mathbf{w}^f,\{L_j\}_{j=1}^n,\pi)$ in terms of the total variation distance $D_{\mathrm{TV}}$ and a strictly convex $f$. The proof can be found in Section \ref{subsubsec:pfcentroidconvalgo}.

\begin{theorem}[$\mathcal{O}\left(1/\sqrt{t}\right)$ convergence rate of the weighted information centroid generated by Algorithm \ref{algo:subg}]\label{thm:centroidconvalgo}
	Consider Algorithm \ref{algo:subg} applied to the family $\{L_j\}_{j=1}^n$ satisfying \eqref{eq:Liclass} under a strictly convex $f$ with its output $(\mathbf{w}^{(i)})_{i=1}^t$ and the notations therein, and the stepsize is chosen to be
	$$\eta = \sqrt{\dfrac{n}{tB}}.$$ 
	We have
	\begin{align}
		D_{\mathrm{TV}}(M^f_n(\overline{\mathbf{w}}^t,\{L_j\}_{j=1}^n,\pi)||M^f_n(\mathbf{w}^f,\{L_j\}_{j=1}^n,\pi)) = \mathcal{O}\left(\dfrac{1}{\sqrt{t}}\right),
	\end{align}
	where the constant in front on the right hand side depends on all the parameters except $t$, that is, the constant depends on $(f,\{L_j\}_{j=1}^n,\pi,n,\mathcal{X})$.
\end{theorem}

\section{Proofs of the main results}\label{sec:proofsmain}

\subsection{Proof of Theorem \ref{thm:existuniquecentroid}}\label{subsubsec:pfexistunique}

The proof is a generalization of \citet[Proposition $1.5$]{DM09} and \citet[Theorem $3.9$]{CW23} to the setting of weighted information centroid. Pick an arbitrary total ordering on $\mathcal{X}$ with strict inequality being denoted by $\prec$. For $i \in \llbracket n \rrbracket$, we write
\begin{align*}
	a &= a(x,y) = \pi(x)M(x,y), \quad	a^{\prime} = a^{\prime}(y,x) = \pi(y)M(y,x), \\
	\beta_i &= \beta_i(x,y) = \pi(x) L_i(x,y), \quad	\beta_i^{\prime} = \beta_i^{\prime}(y,x) = \pi(y)L_i(y,x).
\end{align*}
We note that $M \in \mathcal{L}(\pi)$ gives $a = a^{\prime}$. Using this, we see that
\begin{align*}
	\sum_{i=1}^n w_i D_f(M || L_i) &= \sum_{i=1}^n \sum_{x \prec y} w_i \pi(x) L_i(x,y) f\left(\dfrac{M(x,y)}{L_i(x,y)}\right) +  w_i \pi(y) L_i(y,x) f\left(\dfrac{M(y,x)}{L_i(y,x)}\right) \\
	&= \sum_{i=1}^n \sum_{x \prec y} w_i \beta_i f\left(\dfrac{a}{\beta_i}\right) + w_i \beta_i^{\prime} f\left(\dfrac{a}{\beta_i^{\prime}}\right) \\
	&= \sum_{x \prec y} \sum_{i=1}^n  w_i \beta_i f\left(\dfrac{a}{\beta_i}\right) + w_i \beta_i^{\prime} f\left(\dfrac{a}{\beta_i^{\prime}}\right) \\
	&= \sum_{x \prec y} \sum_{\{i;~w_i>0 ~\textrm{and}~\beta_i>0 ~\textrm{or}~\beta_i^{\prime}>0\}}  w_i \beta_i f\left(\dfrac{a}{\beta_i}\right) + w_i \beta_i^{\prime} f\left(\dfrac{a}{\beta_i^{\prime}}\right) \\
	&=: \sum_{x \prec y} \Phi_{\mathbf{w},\beta_1,\ldots,\beta_n,\beta_1^{\prime},\ldots,\beta_n^{\prime}}(a).
\end{align*}
To minimize with respect to $M$, we are led to minimize the summand above $\phi := \Phi_{\mathbf{w},\beta_1,\ldots,\beta_n,\beta_1^{\prime},\ldots,\beta_n^{\prime}} : \mathbb{R}_+ \to \mathbb{R}_+$, where $\mathbf{w} \in \mathcal{S}_n$ and $(\beta_1,\ldots,\beta_n,\beta_1^{\prime},\ldots,\beta_n^{\prime}) \in \mathbb{R}_+^{2n}$ are assumed to be fixed. The summation in the expression of $\phi$ is non-empty in view of the assumptions in the Theorem.

As $\phi$ and $f$ are convex, we denote by $\phi_+^{\prime}$ and $f_+^{\prime}$ to be their right derivative respectively. It suffices to show the existence of $a_* > 0$ such that for all $a \in \mathbb{R}_+$,
\begin{align}\label{eq:alpha*}
	\phi_+^{\prime}(a) = \begin{cases}
		< 0, \quad \textrm{ if } \quad a < a_*, \\
		> 0, \quad \textrm{ if } \quad a > a_*. 
	\end{cases}
\end{align}
Now, we compute that for all $a \in \mathbb{R}_+$,
\begin{align*}
	\phi_+^{\prime}(a) = \sum_{\{i;~\beta_i>0 ~\textrm{and}~\beta_i^{\prime}>0\}}   w_i f^{\prime}_+\left(\dfrac{a}{\beta_i}\right) + w_i f^{\prime}_+\left(\dfrac{a}{\beta_i^{\prime}}\right) + \sum_{\{i;~\beta_i>0 ~\textrm{and}~\beta_i^{\prime}=0\}}   w_i f^{\prime}_+\left(\dfrac{a}{\beta_i}\right)  + \sum_{\{i;~\beta_i=0 ~\textrm{and}~\beta_i^{\prime}>0\}}   w_i f^{\prime}_+\left(\dfrac{a}{\beta_i^{\prime}}\right).
\end{align*}
As $f^{\prime}(1) = 0$ and $f$ is strictly convex, for sufficiently small $a > 0$ $\phi_+^{\prime}(a) < 0$ while for sufficiently large $a > 0$ $\phi_+^{\prime}(a) > 0$ and $\phi_+^{\prime}$ is increasing, we conclude that there exists a unique $a_* > 0$ such that \eqref{eq:alpha*} is satisfied.

If we replace $f$ by $f^*$ in the above proof, and by noting that $f^*$ is also a strictly convex function with $f^*(1) = f^{*\prime}(1) = 0$, the existence and uniqueness of $M^{f^*}_n$ is shown.

\subsection{Proof of Theorem \ref{thm:examplecentroid}}\label{subsubsec:examplecentroid}

We shall only prove \eqref{eq:cephellinger} as the rest follows exactly the same computation procedure with different choices of $f$. Pick an arbitrary total ordering on $\mathcal{X}$ with strict inequality being denoted by $\prec$. For $i = 1,\ldots,n$, we also write
\begin{align*}
	a &= a(x,y) = \pi(x)M(x,y), \quad	a^{\prime} = a^{\prime}(y,x) = \pi(y)M(y,x), \\
	\beta_i &= \beta_i(x,y) = \pi(x) L_i(x,y), \quad	\beta_i^{\prime} = \beta_i^{\prime}(y,x) = \pi(y)L_i(y,x).
\end{align*}
The $\pi$-reversibility of $M$ yields $a = a^{\prime}$, which leads to
\begin{align*}
	\sum_{i=1}^n w_i D_f(M || L_i) &= \sum_{i=1}^n \sum_{x \prec y} w_i\left(\pi(x) L_i(x,y) f\left(\dfrac{M(x,y)}{L_i(x,y)}\right) +  \pi(y) L_i(y,x) f\left(\dfrac{M(y,x)}{L_i(y,x)}\right)\right) \\
	&= \sum_{i=1}^n \sum_{x \prec y} w_i\left(a - 2 \sqrt{a \beta_i} + \beta_i + a^{\prime} - 2 \sqrt{a^{\prime} \beta_i^{\prime}} + \beta_i^{\prime}\right) \\
	&= \sum_{x \prec y} \sum_{i=1}^n  w_i\left(2a - 2 \sqrt{a \beta_i}  - 2 \sqrt{a \beta_i^{\prime}} + \beta_i + \beta_i^{\prime}\right).
\end{align*}
Next, we aim to minimize each term in the sum, leading us to minimize the strictly convex function of $a$:
$$a \mapsto \sum_{i=1}^n  w_i \left(2a - 2 \sqrt{a \beta_i}  - 2 \sqrt{a \beta_i^{\prime}}\right).$$
Differentiating the above expression with respect to $a$ yields
\begin{align*}
	M^{f}_n(x,y) =\left(\sum_{i=1}^n w_i \sqrt{M^{f}(L_i,\pi)(x,y)}\right)^2.
\end{align*}

\subsection{Proof of Theorem \ref{thm:mainpure}}\label{subsec:pfmainpure}

We first prove item \eqref{it:sd} and \eqref{it:ccic}. To see that the strong duality holds for \eqref{problem:minimax2} and the dual optimal is attained, we shall show that the Slater's constraint qualification (\cite[Theorem $A.1$]{Beck2017} or \cite[Section $5.2.3$]{Boyd2004}) is verified, which amounts to \textcolor{black}{proving} that the constraints in \eqref{problem:minimax2} are strictly feasible. We take
\begin{align*}
	M &= M^f(L_1,\pi), \\
	r &= \max_{i \in \llbracket n \rrbracket} D_f(M||L_i) + 1 >  D_f(M||L_l),
\end{align*}
for all $l \in \llbracket n \rrbracket$, and hence the pair $(M,r)$ is strictly feasible. This proves the first part of item \eqref{it:sd}. As the strong duality holds, using the optimality conditions under strong duality \cite[Theorem $A.2$]{Beck2017}, it proves item \eqref{it:ccic}. In particular, the complementary slackness conditions entail that for all $i \in \llbracket n \rrbracket$,
$$w_i (D_f(M||L_i) - r) = 0,$$
which is equivalent to 
\begin{align*}
	D_f(M||L_i)  \begin{cases}
		= r  & \text{ if } w_i > 0,\\
		\leq r & \text{ if } w_i = 0.
	\end{cases}
\end{align*}
We return to show the second part of item \eqref{it:sd}. By the first part and item \eqref{it:ccic}, we have
$$\overline{v}^f = \sum_{i=1}^n w^f_i D_f(M^{f}_n(\mathbf{w}^f,\{L_i\}_{i=1}^n,\pi)||L_i) = D_f(M^{f}_n(\mathbf{w}^f,\{L_i\}_{i=1}^n,\pi)||L_l),$$
where $l \in \llbracket n \rrbracket$ is such that $w^f_l > 0$.

Next, we prove item \eqref{it:cl}. As the Lagrangian dual is concave \cite[Section $5.1.2$]{Boyd2004}, in our case this amounts to
$$\mathbf{w} \mapsto \min_{r \geq 0, M \in \mathcal{L}(\pi)} \mathsf{L}(r,M,\mathbf{w}) = \sum_{i=1}^n w_i D_f(M^f_n||L_i)$$
is concave.

We proceed to prove item \eqref{it:uniquecc}. Using item \eqref{it:sd}, we see that the pair $(M^{f}_n(\mathbf{w}_1,\{L_i\}_{i=1}^n,\pi),\overline{v}^f)$ minimizes the primal problem and $\mathbf{w}_2 = (w_{2,i})_{i=1}^n$ maximizes the dual problem, and hence by item \eqref{it:ccic}, this gives
$$M^{f}_n(\mathbf{w}_1,\{L_i\}_{i=1}^n,\pi) \in \argmin_{M \in \mathcal{L}(\pi)} \mathsf{L}(r,M,\mathbf{w}_2) = \argmin_{M \in \mathcal{L}(\pi)} \sum_{i=1}^n w_{2,i} D_f(M||L_i).$$
As $f$ is strictly convex, the minimization problem on the right hand side admits a unique minimizer $M^{f}_n(\mathbf{w}_2,\{L_i\}_{i=1}^n,\pi)$ (Theorem \ref{thm:existuniquecentroid}), and hence
$M^{f}_n(\mathbf{w}_1,\{L_i\}_{i=1}^n,\pi) = M^{f}_n(\mathbf{w}_2,\{L_i\}_{i=1}^n,\pi)$.

Next, we prove item \eqref{it:vfequals}. We first show sufficiency. Using item \eqref{it:ccic}, we note that
$$\overline{v}^f = D_f(M^{f}_n(\mathbf{w}^f,\{L_i\}_{i=1}^n,\pi)||L_j) = \min_{M \in \mathcal{L}(\pi)} D_f(M||L_j) \leq \max_{l \in \llbracket n \rrbracket} \min_{M \in \mathcal{L}(\pi)} D_f(M||L_l) = \underline{v}^f.$$
This gives $\overline{v}^f = \underline{v}^f$ since $\overline{v}^f \geq \underline{v}^f$. To show necessity, let 
$$j = \argmax_{l \in \llbracket n \rrbracket} \min_{M \in \mathcal{L}(\pi)} D_f(M||L_l),$$
and we take $\mathbf{w}^f$ to be the standard unit vector at the $j$-th coordinate, that is, $w^f_j = 1$ and $0$ otherwise. Since $\overline{v}^f = \underline{v}^f$, this choice of $\mathbf{w}^f$ maximizes the dual problem \eqref{problem:minimax2dual} and 
\begin{align*}
	M^{f}_n(\mathbf{w}^f,\{L_i\}_{i=1}^n,\pi) \in \argmin_{M \in \mathcal{L}(\pi)} D_f(M||L_j)
\end{align*}
clearly holds.

Finally, we prove item \eqref{it:saddlept}. As the saddle point property holds we have
$$\overline{v}^f = \max_{i \in \llbracket n \rrbracket} D_f(M^f(L_i,\pi),L_i).$$
The stated pair in the theorem statement is thus a saddle point.

\subsection{Proof of Proposition \ref{prop:reducefinite}}\label{subsubsec:pfreducefinite}

First, as $D_f(\cdot||\cdot)$ is jointly convex in its arguments, and maximization of a convex function over a convex hull occur at an extreme point, we see that
$$\sup_{L \in \mathrm{conv}((L_i)_{i=1}^n)} D_f(M||L) = \max_{i \in \llbracket n \rrbracket} D_f(M||L_i).$$
Taking inf on both sides over $M \in \mathcal{L}(\pi)$ gives the desired result for $\overline{v}^f$. 

For $\underline{v}^f$, we note that the mapping
$$L \mapsto \inf_{M \in \mathcal{L}(\pi)} D_f(M||L)$$
is convex as $D_f(\cdot||\cdot)$ is jointly convex and partial minimization of convex function preserves convexity \cite[Section $3.2.5$]{Boyd2004}. Using again the property that maximization of a convex function over a convex hull occur at an extreme point, we have
$$\underline{v}^f(\mathrm{conv}((L_i)_{i=1}^n),\pi) = \sup_{L \in \mathrm{conv}((L_i)_{i=1}^n)} \inf_{M \in \mathcal{L}(\pi)}  D_f(M||L) = \max_{i \in \llbracket n \rrbracket} \inf_{M \in \mathcal{L}(\pi)}  D_f(M||L_i) = \underline{v}^f(\{L_i\}_{i=1}^n,\pi).$$

\subsection{Proof of Corollary \ref{cor:pureNash}}\label{subsubsec:pfpureNash}

For the first item, it is simply a restatement of the saddle point result Corollary \ref{cor:convexhull} in the context of the two-person game and Nash equilibrium.

For the second item, as the index $l$ is unique, using the fact that $f$ is strictly convex $M^f(L_l,\pi)$ is unique, and we thus have a unique saddle point and hence Nash equilibrium. For the other direction, suppose that $l_1, l_2 \in  \argmax_{i \in \llbracket n \rrbracket} D_f(M^f(L_i,\pi),L_i)$, then both $(M^f(L_{l_1},\pi),L_{l_1})$ and $(M^f(L_{l_2},\pi),L_{l_2})$ are saddle points or Nash equilibria.

\subsection{Proof of Proposition \ref{prop:minimaxeqmixed}}\label{subsubsec:pfminimaxineq}

We first prove the inequalities in \eqref{eq:minimaxeqmixed1}. The inequality $\overline{V}^f \geq \underline{V}^f$ is obvious. To see the first inequality, we note that
\begin{align*}
	\int_{\mathcal{B}} D_f(M||L)\, \mu(dL) \leq \sup_{L \in \mathcal{B}} D_f(M||L).
\end{align*}
The desired result follows by taking $\sup_{\mu \in \mathcal{P}(\mathcal{B})}$ then $\inf_{M \in \mathcal{L}(\pi)}$.

Next, we prove the inequalities in \eqref{eq:minimaxeqmixed2}, and it only remains to show the rightmost inequality. We observe that
\begin{align*}
	\underline{V}^f(\mathcal{B},\pi) &= \sup_{\mu \in \mathcal{P}(\mathcal{B})} \inf_{M \in \mathcal{L}(\pi)} \int_{\mathcal{B}} D_f(M||L)\, \mu(dL) \\
	&\geq \max_{\mathbf{w} \in \mathcal{S}_n} \inf_{M \in \mathcal{L}(\pi)} \sum_{i=1}^n w_i D_f(M||L_i)  \\
	&= \max_{\mathbf{w} \in \mathcal{S}_n} \sum_{i=1}^n w_i D_f(M^f_n||L_i),
\end{align*}
which completes the proof.

\subsection{Proof of Theorem \ref{thm:subgh}}\label{subsubsec:pfsubgh}

First, by the definition of weighted information centroid, we have
$$\sum_{i=1}^n w_i D_f(M^{f}_n(\mathbf{w},\{L_j\}_{j=1}^n,\pi)||L_i) \leq \sum_{i=1}^n w_i D_f(M^{f}_n(\mathbf{v},\{L_j\}_{j=1}^n,\pi)||L_i).$$
Multiplying both sides by $-1$ followed by substracting $h(\mathbf{v})$, we obtain
\begin{align*}
	h(\mathbf{w}) - h(\mathbf{v}) &= - \sum_{i=1}^n w_i D_f(M^{f}_n(\mathbf{w},\{L_j\}_{j=1}^n,\pi)||L_i) + \sum_{i=1}^n v_i D_f(M^{f}_n(\mathbf{v},\{L_j\}_{j=1}^n,\pi)||L_i) \\
	&\geq \sum_{i=1}^n -D_f(M^{f}_n(\mathbf{v},\{L_j\}_{j=1}^n,\pi)||L_i)(w_i-v_i) \\
	&= \sum_{i=1}^n -D_f(M^{f}_n(\mathbf{v},\{L_j\}_{j=1}^n,\pi)||L_i)(w_i-v_i) + \sum_{i=1}^n D_f(M^{f}_n(\mathbf{v},\{L_j\}_{j=1}^n,\pi)||L_n)(w_i-v_i) \\
	&= \sum_{i=1}^n g_i (w_i - v_i),
\end{align*}
where the second equality follows from $\mathbf{w}, \mathbf{v} \in \mathcal{S}_n$.

Next, we proceed to prove \eqref{eq:subgnormbd}. First we note that for any $i \in \llbracket n \rrbracket$, we have
\begin{align*}
	D_f(M^{f}_n(\mathbf{w},\{L_j\}_{j=1}^n,\pi)||L_i) &= \sum_{x \in \mathcal{X}} \pi(x) \sum_{y \in \mathcal{X}\backslash\{x\}} L_i(x,y) f\left(\dfrac{M^{f}_n(\mathbf{w},\{L_j\}_{j=1}^n,\pi)(x,y)}{L_i(x,y)}\right) \\
	&\leq |\mathcal{X}| \sup_{\mathbf{v} \in \mathcal{S}_n;~i \in \llbracket n \rrbracket;~L_i(x,y)>0} L_i(x,y) f\left(\dfrac{M^f_n(\mathbf{v},\{L_j\}_{j=1}^n,\pi)(x,y)}{L_i(x,y)}\right),
\end{align*}
and hence
\begin{align*}
	\norm{\mathbf{g}}_2^2 = \sum_{i=1}^n g_i^2 &\leq \sum_{i=1}^n \max\bigg\{D_f(M^{f}_n(\mathbf{w},\{L_j\}_{j=1}^n,\pi)||L_n)^2, D_f(M^{f}_n(\mathbf{w},\{L_j\}_{j=1}^n,\pi)||L_i)^2\bigg\} \\
	&\leq n \max_{l \in \llbracket n \rrbracket} D_f(M^f_n(\mathbf{v},\{L_j\}_{j=1}^n,\pi)||L_l)^2 \leq B.
\end{align*}

\subsection{Proof of Theorem \ref{thm:subgconv}}\label{subsubsec:pfsubgconv}

We first prove \eqref{eq:subgconvbd}. For $i \in \llbracket t \rrbracket$, we have

\begin{align*}
	\norm{\mathbf{w}^{(i+1)} - \mathbf{w}^f}_2^2 &\leq \norm{\mathbf{v}^{(i+1)} - \mathbf{w}^f}_2^2\\
											   &= \norm{\mathbf{w}^{(i)} - \mathbf{w}^f - \eta \mathbf{g}(\mathbf{w}^{(i)})}_2^2 \\
											   &= \norm{\mathbf{w}^{(i)} - \mathbf{w}^f}_2^2 + \eta^2 \norm{\mathbf{g}(\mathbf{w}^{(i)})}_2^2 - 2 \eta \mathbf{g}(\mathbf{w}^{(i)}) \cdot (\mathbf{w}^{(i)} - \mathbf{w}^f) \\
											   &\leq \norm{\mathbf{w}^{(i)} - \mathbf{w}^f}_2^2 + \eta^2 B - 2 \eta \mathbf{g}(\mathbf{w}^{(i)}) \cdot (\mathbf{w}^{(i)} - \mathbf{w}^f),
\end{align*}
where the first inequality follows from the fact that $\mathbf{w}^{(i+1)}$ is the projection of $\mathbf{v}^{(i+1)}$ onto the simplex $\mathcal{S}_n$, while the second inequality follows from the definition of $B$ in \eqref{eq:subgnormbd}. Upon rearranging and using the definition of subgradient, this leads to
\begin{align*}
	h(\mathbf{w}^{(i)}) - h(\mathbf{w}^f) &\leq \mathbf{g}(\mathbf{w}^{(i)}) \cdot (\mathbf{w}^{(i)} - \mathbf{w}^f)\\
	&= \dfrac{1}{2 \eta} \left(\norm{\mathbf{w}^{(i)} - \mathbf{w}^f}_2^2 - \norm{\mathbf{w}^{(i+1)} - \mathbf{w}^f}_2^2\right)+\dfrac{\eta B}{2}.
\end{align*}
Now, we sum $i$ from $1$ to $t$ to give
\begin{align*}
	\sum_{i=1}^t h(\mathbf{w}^{(i)}) - h(\mathbf{w}^f) &\leq \dfrac{1}{2 \eta} \left(\norm{\mathbf{w}^{(1)} - \mathbf{w}^f}_2^2 - \norm{\mathbf{w}^{(t+1)} - \mathbf{w}^f}_2^2\right)+\dfrac{\eta B t}{2} \\
	&\leq \dfrac{1}{2 \eta} \norm{\mathbf{w}^{(1)} - \mathbf{w}^f}_2^2 +\dfrac{\eta B t}{2} \leq \dfrac{1}{2 \eta} n +\dfrac{\eta B t}{2},
\end{align*}
where in the last inequality we use the fact that $\mathbf{w}^{(1)},\mathbf{w}^f \in \mathcal{S}_n$. Dividing both sides by $t$ followed by the convexity of $h$ (see \eqref{def:h}) yields
\begin{align*}
	h\left(\dfrac{1}{t}\sum_{i=1}^t \mathbf{w}^{(i)}\right) - h\left(\mathbf{w}^{f}\right) \leq \dfrac{1}{t} \left(\sum_{i=1}^t h(\mathbf{w}^{(i)}) - h(\mathbf{w}^f)\right) \leq \dfrac{n}{2\eta t} + \dfrac{\eta B}{2}.
\end{align*}
The right hand side as a function of the stepsize $\eta$ is minimized when we take $\eta = \sqrt{\dfrac{n}{tB}}$.

\subsection{Proof of Theorem \ref{thm:centroidconvalgo}}\label{subsubsec:pfcentroidconvalgo}

We first state the following lemma in the setting of this Theorem, which is of independent interest. It can be considered as an extension of \cite[equation $(3.25)$]{C72} to finite Markov chains:

\begin{lemma}\label{lem:important}
	There exists a constant $\infty > C = C(f,\{L_i\}_{i=1}^n,\pi,n,\mathcal{X}) > 0$ such that
	\begin{align*}
		D_{\mathrm{TV}}&(M^f_n(\overline{\mathbf{w}}^t,\{L_i\}_{i=1}^n,\pi)||M^f_n(\mathbf{w}^f,\{L_i\}_{i=1}^n,\pi)) \\
		&\leq C \left(\sum_{i=1}^n w^f_i D_f(M^{f}_n(\mathbf{w}^f,\{L_i\}_{i=1}^n,\pi)||L_i) - \sum_{i=1}^n \overline{w}^t_i D_f(M^{f}_n(\overline{\mathbf{w}}^t,\{L_i\}_{i=1}^n,\pi)||L_i)\right) \\
		&= C \left(h\left(\overline{\mathbf{w}}^t\right) - h\left(\mathbf{w}^{f}\right)\right).
	\end{align*}
\end{lemma}

Once we have Lemma \ref{lem:important}, the desired result is obtained, since from Theorem \ref{thm:subgconv} with the same choice of stepsize in Algorithm \ref{algo:subg} we know that
$$D_{\mathrm{TV}}(M^f_n(\overline{\mathbf{w}}^t,\{L_i\}_{i=1}^n,\pi)||M^f_n(\mathbf{w}^f,\{L_i\}_{i=1}^n,\pi)) \leq C \left(h\left(\overline{\mathbf{w}}^t\right) - h\left(\mathbf{w}^{f}\right)\right) = \mathcal{O}\left(\dfrac{1}{\sqrt{t}}\right).$$

It thus remains to prove Lemma \ref{lem:important}. In the remaining of this proof, to minimize notational overhead we write
$$M^f_n(\overline{\mathbf{w}}^t) = M^f_n(\overline{\mathbf{w}}^t,\{L_i\}_{i=1}^n,\pi), \quad M^f_n(\mathbf{w}^f) = M^f_n(\mathbf{w}^f,\{L_i\}_{i=1}^n,\pi).$$
Using the strict convexity of $f$ and \eqref{eq:Liclass}, according to \cite[equation $(3.25)$]{C72} we have that
\begin{align*}
	D_{\mathrm{TV}}&(M^f_n(\overline{\mathbf{w}}^t)||M^f_n(\mathbf{w}^f)) \\
	&\leq C \left(\sum_{i=1}^n \overline{w}^t_i D_f(M^{f}_n(\mathbf{w}^f,\{L_i\}_{i=1}^n,\pi)||L_i) - \sum_{i=1}^n \overline{w}^t_i D_f(M^{f}_n(\overline{\mathbf{w}}^t,\{L_i\}_{i=1}^n,\pi)||L_i)\right) \\
	&\leq C \left(\max_{i \in \llbracket n \rrbracket} D_f(M^{f}_n(\mathbf{w}^f,\{L_i\}_{i=1}^n,\pi)||L_i) + h\left(\overline{\mathbf{w}}^t\right)\right) \\
	&= C \left(- h\left(\mathbf{w}^{f}\right) + h\left(\overline{\mathbf{w}}^t\right)\right),
\end{align*}
where the last equality follows from the complementary slackness in Theorem \ref{thm:mainpure}.

\section*{Acknowledgements}

\textcolor{black}{The authors gratefully acknowledge the careful reading and constructive comments of a reviewer that has improved the quality and presentation of the manuscript. Michael Choi would like to thank Laurent Miclo for raising a question about possible game-theoretic interpretations of $P_{-\infty}$ and $P_{\infty}$ when he visited Toulouse School of Economics.} For a possible answer to this question, the readers are referred to Example \ref{ex:1mix}. \textcolor{black}{Michael Choi acknowledges the financial support of the projects A-8001061-00-00, NUSREC-HPC-00001, NUSREC-CLD-00001, A-0000178-01-00, A-0000178-02-00 and A-8003574-00-00 at National University of Singapore.} Part of this research was supported when Geoffrey Wolfer was part of the Special Postdoctoral Researcher Program (SPDR) of RIKEN and by the Japan Society for the Promotion of Science KAKENHI under Grant 23K13024.
\bibliographystyle{abbrvnat}
\bibliography{thesis}

\begin{thebibliography}{33}
\providecommand{\natexlab}[1]{#1}
\providecommand{\url}[1]{\texttt{#1}}
\expandafter\ifx\csname urlstyle\endcsname\relax
  \providecommand{\doi}[1]{doi: #1}\else
  \providecommand{\doi}{doi: \begingroup \urlstyle{rm}\Url}\fi

\bibitem[Beck(2017)]{Beck2017}
A.~Beck.
\newblock \emph{First-order methods in optimization}, volume~25 of
  \emph{MOS-SIAM Series on Optimization}.
\newblock Society for Industrial and Applied Mathematics (SIAM), Philadelphia,
  PA; Mathematical Optimization Society, Philadelphia, PA, 2017.

\bibitem[Bierkens(2016)]{Bie16}
J.~Bierkens.
\newblock Non-reversible {M}etropolis-{H}astings.
\newblock \emph{Stat. Comput.}, 26\penalty0 (6):\penalty0 1213--1228, 2016.

\bibitem[Boyd and Vandenberghe(2004)]{Boyd2004}
S.~Boyd and L.~Vandenberghe.
\newblock \emph{Convex optimization}.
\newblock Cambridge University Press, Cambridge, 2004.

\bibitem[Candan(2020)]{C20}
C.~Candan.
\newblock Chebyshev center computation on probability simplex with
  $\alpha$-divergence measure.
\newblock \emph{IEEE Signal Processing Letters}, 27:\penalty0 1515--1519, 2020.

\bibitem[Choi(2020)]{Choi16}
M.~C. Choi.
\newblock Metropolis-{H}astings reversiblizations of non-reversible {M}arkov
  chains.
\newblock \emph{Stochastic Processes and their Applications}, 130\penalty0
  (2):\penalty0 1041 -- 1073, 2020.

\bibitem[Choi and Huang(2020)]{CH18}
M.~C. Choi and L.-J. Huang.
\newblock On hitting time, mixing time and geometric interpretations of
  {M}etropolis-{H}astings reversiblizations.
\newblock \emph{J. Theoret. Probab.}, 33\penalty0 (2):\penalty0 1144--1163,
  2020.

\bibitem[Choi and Wolfer(2024)]{CW23}
M.~C.~H. Choi and G.~Wolfer.
\newblock Systematic approaches to generate reversiblizations of {M}arkov
  chains.
\newblock \emph{IEEE Trans. Inform. Theory}, 70\penalty0 (5):\penalty0
  3145--3161, 2024.

\bibitem[Condat(2016)]{Condat16}
L.~Condat.
\newblock Fast projection onto the simplex and the {$l_1$} ball.
\newblock \emph{Math. Program.}, 158\penalty0 (1-2):\penalty0 575--585, 2016.

\bibitem[Csisz\'{a}r(1972)]{C72}
I.~Csisz\'{a}r.
\newblock A class of measures of informativity of observation channels.
\newblock \emph{Period. Math. Hungar.}, 2:\penalty0 191--213, 1972.

\bibitem[Diaconis and Miclo(2009)]{DM09}
P.~Diaconis and L.~Miclo.
\newblock On characterizations of {M}etropolis type algorithms in continuous
  time.
\newblock \emph{ALEA Lat. Am. J. Probab. Math. Stat.}, 6:\penalty0 199--238,
  2009.

\bibitem[Diaconis et~al.(2000)Diaconis, Holmes, and Neal]{DHN00}
P.~Diaconis, S.~Holmes, and R.~M. Neal.
\newblock Analysis of a nonreversible {M}arkov chain sampler.
\newblock \emph{Ann. Appl. Probab.}, 10\penalty0 (3):\penalty0 726--752, 2000.

\bibitem[Eldar et~al.(2008)Eldar, Beck, and Teboulle]{EBT08}
Y.~C. Eldar, A.~Beck, and M.~Teboulle.
\newblock A minimax {C}hebyshev estimator for bounded error estimation.
\newblock \emph{IEEE Trans. Signal Process.}, 56\penalty0 (4):\penalty0
  1388--1397, 2008.

\bibitem[Fill(1991)]{Fill91}
J.~A. Fill.
\newblock Eigenvalue bounds on convergence to stationarity for nonreversible
  {M}arkov chains, with an application to the exclusion process.
\newblock \emph{Ann. Appl. Probab.}, 1\penalty0 (1):\penalty0 62--87, 1991.

\bibitem[Gushchin and Zhdanov(2006)]{GZ06}
A.~A. Gushchin and D.~A. Zhdanov.
\newblock A minimax result for {$f$}-divergences.
\newblock In \emph{From stochastic calculus to mathematical finance}, pages
  287--294. Springer, Berlin, 2006.

\bibitem[Gustafson(1998)]{G98}
P.~Gustafson.
\newblock A guided walk metropolis algorithm.
\newblock \emph{Statistics and computing}, 8\penalty0 (4):\penalty0
  357{\textendash}364, 1998.

\bibitem[Haussler(1997)]{H97}
D.~Haussler.
\newblock A general minimax result for relative entropy.
\newblock \emph{IEEE Trans. Inform. Theory}, 43\penalty0 (4):\penalty0
  1276--1280, 1997.

\bibitem[Haussler and Opper(1997)]{HO97}
D.~Haussler and M.~Opper.
\newblock Mutual information, metric entropy and cumulative relative entropy
  risk.
\newblock \emph{Ann. Statist.}, 25\penalty0 (6):\penalty0 2451--2492, 1997.

\bibitem[Horn and Johnson(2013)]{HJ13}
R.~A. Horn and C.~R. Johnson.
\newblock \emph{Matrix analysis}.
\newblock Cambridge University Press, Cambridge, second edition, 2013.

\bibitem[Hwang et~al.(1993)Hwang, Hwang-Ma, and Sheu]{Hwang93}
C.-R. Hwang, S.-Y. Hwang-Ma, and S.~J. Sheu.
\newblock Accelerating {G}aussian diffusions.
\newblock \emph{Ann. Appl. Probab.}, 3\penalty0 (3):\penalty0 897--913, 1993.

\bibitem[Hwang et~al.(2005)Hwang, Hwang-Ma, and Sheu]{Hwang05}
C.-R. Hwang, S.-Y. Hwang-Ma, and S.-J. Sheu.
\newblock Accelerating diffusions.
\newblock \emph{Ann. Appl. Probab.}, 15\penalty0 (2):\penalty0 1433--1444,
  2005.

\bibitem[Jansen and Kurt(2014)]{JK14}
S.~Jansen and N.~Kurt.
\newblock On the notion(s) of duality for {M}arkov processes.
\newblock \emph{Probab. Surv.}, 11:\penalty0 59--120, 2014.

\bibitem[Kamatani and Song(2023)]{KS23}
K.~Kamatani and X.~Song.
\newblock Non-reversible guided {M}etropolis kernel.
\newblock \emph{J. Appl. Probab.}, 60\penalty0 (3):\penalty0 955--981, 2023.

\bibitem[Karlin and Peres(2017)]{KarlinPeres17}
A.~R. Karlin and Y.~Peres.
\newblock \emph{Game theory, alive}.
\newblock American Mathematical Society, Providence, RI, 2017.

\bibitem[Keilson(1979)]{K79}
J.~Keilson.
\newblock \emph{Markov chain models---rarity and exponentiality}, volume~28 of
  \emph{Applied Mathematical Sciences}.
\newblock Springer-Verlag, New York-Berlin, 1979.

\bibitem[Maschler et~al.(2020)Maschler, Solan, and Zamir]{MSZ20}
M.~Maschler, E.~Solan, and S.~Zamir.
\newblock \emph{Game theory}.
\newblock Cambridge University Press, Cambridge, second edition, 2020.

\bibitem[Miclo(1997)]{M97}
L.~Miclo.
\newblock Remarques sur l'hypercontractivit\'{e} et l'\'{e}volution de
  l'entropie pour des cha\^{i}nes de {M}arkov finies.
\newblock In \emph{S\'{e}minaire de {P}robabilit\'{e}s, {XXXI}}, volume 1655 of
  \emph{Lecture Notes in Math.}, pages 136--167. Springer, Berlin, 1997.

\bibitem[Paulin(2015)]{Paulin15}
D.~Paulin.
\newblock Concentration inequalities for {M}arkov chains by {M}arton couplings
  and spectral methods.
\newblock \emph{Electron. J. Probab.}, 20:\penalty0 no. 79, 1--32, 2015.

\bibitem[Rey-Bellet and Spiliopoulos(2016)]{BS16}
L.~Rey-Bellet and K.~Spiliopoulos.
\newblock Improving the convergence of reversible samplers.
\newblock \emph{J. Stat. Phys.}, 164\penalty0 (3):\penalty0 472--494, 2016.

\bibitem[Rosenthal and Rosenthal(2015)]{RR15}
J.~S. Rosenthal and P.~Rosenthal.
\newblock Spectral bounds for certain two-factor non-reversible {MCMC}
  algorithms.
\newblock \emph{Electron. Commun. Probab.}, 20:\penalty0 no. 91, 10, 2015.

\bibitem[Schr{\"o}dinger(1931)]{S1931}
E.~Schr{\"o}dinger.
\newblock \"uber die umkehrung der naturgesetze.
\newblock \emph{Sitzungsberichte der Preussischen Akademie der Wissenschaften,
  Physikalisch-Mathematische Klasse}, 8\penalty0 (9):\penalty0 144--153, 1931.

\bibitem[van Dijk et~al.(2018)van Dijk, van Brummelen, and Boucherie]{V18}
N.~M. van Dijk, S.~P.~J. van Brummelen, and R.~J. Boucherie.
\newblock Uniformization: Basics, extensions and applications.
\newblock \emph{Perform. Evaluation}, 118:\penalty0 8--32, 2018.

\bibitem[van Erven and Harremo\"{e}s(2014)]{EH14}
T.~van Erven and P.~Harremo\"{e}s.
\newblock R\'{e}nyi divergence and {K}ullback-{L}eibler divergence.
\newblock \emph{IEEE Trans. Inform. Theory}, 60\penalty0 (7):\penalty0
  3797--3820, 2014.

\bibitem[Wolfer and Watanabe(2021)]{WW21}
G.~Wolfer and S.~Watanabe.
\newblock Information geometry of reversible {M}arkov chains.
\newblock \emph{Inf. Geom.}, 4\penalty0 (2):\penalty0 393--433, 2021.

\end{thebibliography}

\end{document}